\documentclass[sn_mathphys,Numbered]{sn_jnl}

\usepackage{graphicx}%
\usepackage{multirow}%
\usepackage{amsmath,amssymb,amsfonts}%
\usepackage{amsthm}%
\usepackage{mathrsfs}%
\usepackage[title]{appendix}%
\usepackage{xcolor}%
\usepackage{textcomp}%
\usepackage{manyfoot}%
\usepackage{booktabs}%
\usepackage{algorithm}%
\usepackage{algorithmicx}%
\usepackage{algpseudocode}%
\usepackage{listings}%
\usepackage{hyperref}

\usepackage[caption=false]{subfig}

\usepackage{tikz}
\usetikzlibrary{positioning}
\usepackage{wrapfig}

\theoremstyle{thmstyleone}%
\theoremstyle{thmstyletwo}%

\theoremstyle{thmstylethree}%

\begin{document}

\title[Article Title]{Ensemble Kalman Filtering for Glacier Modeling}



\abstract{Working with a two-stage ice sheet model, we explore how statistical data assimilation methods can be used to improve predictions of glacier melt and relatedly, sea level rise. We find that the EnKF improves model runs initialized using incorrect initial conditions or parameters, providing us with better models of future glacier melt. We explore the necessary number of observations needed to produce an accurate model run. Further, we determine that the deviations from the truth in output that stem from having few data points in the pre-satellite era can be corrected with modern observation data. Finally, using data derived from our improved model we calculate sea level rise and model storm surges to understand the affect caused by sea level rise.}

\author[]{Emily Corcoran , Logan Knudsen, Talea Mayo, Hannah Park-Kaufmann, Alexander Robel}




\maketitle

\section{Introduction}\label{sec2} 
\indent Research has shown that climate change will likely impact storm surge and make storm surges more severe (\cite{camelo2020projected}). Storm surge occurs when high winds and low pressure from a tropical storm force ocean water into coastal regions. Inundation caused by storm surge has the potential to cause a great deal of damage in coastal regions when a storm makes landfall. Models have been developed in order to help predict storm surge during storm events, and researchers have begun to include future sea-level rise as a factor in these models. When studying the impact of recent tropical storms, \cite{camelo2020projected} found that for all 14 storms simulated in the Gulf of Mexico and Atlantic Ocean, storm inundation increased by an average of 25\%, likely due to the impacts of climate change. This illustrates the importance of modeling sea level rise and therefore improving glacier models by incorporating data assimilation. Doing so increases our understanding of how glacier melt will contribute to sea level rise, and in turn affect storm surge. \\
\indent Sea level rise caused by climate change plays a significant part in this impact. To better model sea-level rise, we turn to glacier modeling, specifically marine-terminating glaciers. These have a natural flow towards the ocean, which contributes to sea level rise (\cite{robel2019marine}). By the year 2300, the Antarctic ice sheet is projected to cause up to 3 meters of sea level rise globally (\cite{robel2015long}). Due to the severe impacts of glacial melting, modeling changes in ice sheets is an important task. There are challenges to modeling sea level rise, as ice sheet instability leads to significant sea-level rise uncertainty (\cite{robel2019marine}).\\
\indent Data assimilation is a method to move models closer to reality using real-world observations by readjusting the model state at specified times (\cite{dataassim}). Data assimilation was initially developed for use in weather models in the late 1990s in order to improve weather forecast accuracy on short-range weather models (\cite{barker2012weather}). These models would typically have a run time of 1-6 hours, and during that process real-time observations would be pulled for this purpose in order to readjust the current model run. Since its introduction, data assimilation has been incorporated into other geoscience models and radars, to name a few uses. \\
\indent Data assimilation can also play a factor in decision making when it comes to collecting field data. In the field of glaciology, data is largely collected via satellite or in-person field data collection. While both methods collect valuable data for modeling glaciers, both are expensive and time-consuming. Using data assimilation can help to inform the glacier modelers and glaciologists who collect data about how to collect data in an efficient way. This can help researchers to more efficiently utilize funding and avoid unnecessarily expensive data collection that does not significantly improve glacier models.

\section{Methods} 
\vspace{-1em} 
\subsection{Glacier Model}
\begin{wrapfigure}{r}{0.6\textwidth}
    \includegraphics[width=0.9\linewidth]{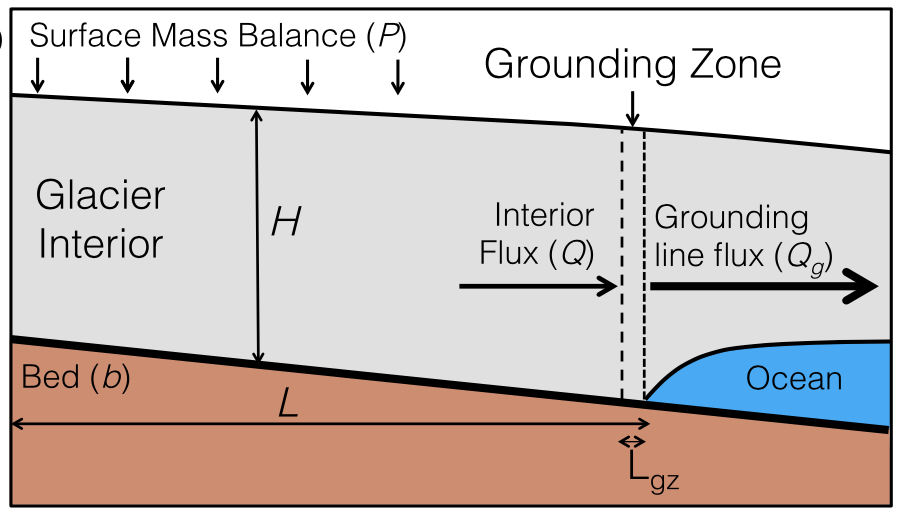}
    \caption{Diagram of a marine-terminating glacier from \cite{robel2018response}.}
    \label{fig:glacierdiagram}
\end{wrapfigure}

\indent We have chosen a simplified two-stage ice sheet model for our exploration (Figure \ref{fig:glacierdiagram}). This model describes the changes in ice mass of marine-terminating glaciers, which may be impacted over time by climate change (\cite{robel2018response}).

A glacier can be represented with a simplified box model that has a length $L$, precipitation $P$, and height and flux at the grounding line, $h_g$ and $Q_g$, respectively. The two-stage model used in this paper incorporates a nested box into the system to more accurately represent a glacier. This new box has a thickness, $H$, and an interior flux, $Q$. The change in length and height of the glacier can be described with these differential equations:
\begin{equation}\label{dH}
  \dfrac{dH}{dt}=P-\dfrac{Q_g}{L}-\dfrac{H}{h_gL}(Q-Q_g)  
\end{equation}
\begin{equation} 
\dfrac{dL}{dt}=\dfrac{1}{h_g}(Q-Q_g). \end{equation}
The following equations are used to calculate the variables used in this system:
\begin{equation}h_g = -\lambda b(L)\end{equation}
\begin{equation}\label{q}Q = \gamma \dfrac{H^{2n+1}}{L^n}\end{equation}
\begin{equation}\label{qg}Q_g = \Omega h_g^\beta,\end{equation}
where $\gamma$ and $\Omega$ are assumed to be constants for our purposes, and $\lambda = \dfrac{\rho_w}{\rho_i}$, where $\rho_w$ is seawater density and $\rho_i$ is glacial ice density.

\begin{wrapfigure}{lH}{0.45\textwidth}
\begin{tabular}{|p{2cm}||p{2cm}|  }
 \hline
 \multicolumn{2}{|c|}{Parameters and Initial Conditions} \\
 \hline
 Parameter & Value\\
 \hline
$smb_o$ & 0.3  \\
$smb_1$ & 0.15 \\
$smb_f$ & 0.0 \\
$H_o$ & 2.18 \\
$L_o$  & 4.44 \\
$b_x$  & -0.001 \\
$sill_{min}$ & 415 \\
$sill_{max}$ & 425 \\
$sill_{slope}$ & 0.01 \\
 \hline
\end{tabular}
\caption{}
\label{fig:initcondtrue}
\end{wrapfigure}

Despite the simplicity of this model, it offers a sufficient approximation of glacier melt using $Q$ and $Q_g$. As a result, findings we learn from this simplified model are worth evaluating in a more complex model. Further, for most of this study we will assume the initial conditions and parameters are what we see in Figure \ref{fig:initcondtrue}. We let $smb_o$, $smb_1$, and $smb_f$ describe the surface mass balance (P in Equation \ref{dH}) at three points over time. $H_o$ and $L_o$ are initial conditions for height and length at year 0. The parameter $b_x$ is the slope of the Earth underneath the glacier. Finally, $sill_{min}$ and $sill_{max}$ describe the start point of the glacial sill (region of reverse slope) with $sill_{slope}$ representing the slope of the sill. For a more in-depth explanation and justification of the model, see \cite{robel2018response}. 

\begin{figure}[H]
    \centering
    \scalebox{0.8}{\begin{tikzpicture}[
    squarednode/.style={rectangle, draw=blue!60, fill=blue!5, very thick,
    minimum size=5mm},
    squarednode2/.style={rectangle, draw=green!60, fill=green!5, very thick, minimum size=5mm},
    squarednode3/.style={rectangle, draw=violet!60, fill=purple!5, very thick, minimum size=5mm},
    ]
    \node[squarednode]      (maintopic)                     {Model};
    \node[squarednode2]         (uppercircle)       [above=of maintopic] {Initial Conditions};
    \node[squarednode3]      (rightsquare)       [left=of maintopic] {Initialize DA Function};
    \node[squarednode3]        (lowercircle)       [below=of maintopic] {Analysis};
    \node[squarednode2]        (lowercircle2)       [right= of maintopic] {Forecast at time $t$};
    
    \draw[-latex] (uppercircle.south) -- (maintopic.north);
    \draw[-latex] (maintopic.west) -- (rightsquare.east);
    \draw[-latex] (rightsquare.south) -- (lowercircle.west);
    \draw[-latex] (maintopic.east) -- (lowercircle2.west);
    \draw[-latex] (lowercircle.north) -- (maintopic.south);
    \end{tikzpicture}}
    \caption{Code Structure Diagram}
    \label{fig:codeDiagram}
\end{figure}
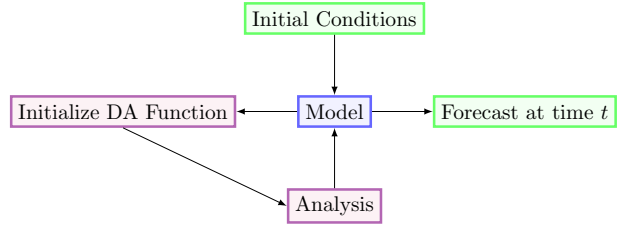

\indent The program used to model the glacier behavior and assimilate the data requires choosing a set of initial conditions (Figure ~\ref{fig:codeDiagram}). Once the initial conditions are input to the model, a Runge-Kutta 4th order method is used to advance the model in time, and the model output can be used along with a data assimilation method. Then the analyzed data from the assimilation is fed back into the model, where time is advanced again. 
    
\indent Next, we give an overview of the data assimilation techniques used in this work. 
\subsection{Kalman Filter}
    \indent The Kalman Filter is a data assimilation technique that uses the model, observations, and corresponding error covariance matrices in order to adjust model output to be closer to reality. The goal of the Kalman Filter is to compute an optimal estimate from a combination of the results of a previous forecast and observations. The key to this is the Kalman Gain, $\textbf{K}_t$, which decides the balance of how much this analysis relies on the model and the observations. 
    
\subsubsection{Ensemble Kalman Filter} 
    \indent The Ensemble Kalman Filter (EnKF) is a nonlinear version of the Kalman filter fit for large problems. The model state is represented by an ensemble of states, and the covariance matrix is replaced by the sample covariance. An analysis is performed for each member of the ensemble. \\
    Consider the following state estimation system:
\begin{equation}
\mathbf{x}^{f,(i)}_{t}=\mathcal{M} \mathbf{x}^{(i)}_{t-1}+w^{(i)}_{t}
\end{equation}
\begin{equation}
\mathbf{P}_{f}= \frac{1}{N-1} \sum_{i=0}^N[\mathbf{x}^{(i)}_t-\mathbf{x}^a_i][\mathbf{x}^{(i)}_t-\mathbf{x}^a_i]'
\end{equation}


Let $\mathcal{M}$ represent our model and $\mathbf{x}^{a,(0)}_t, \mathbf{x}^{a,(1)}_t, \dots, \mathbf{x}^{a,(N)}_t,$ be the ensemble members at time $t$. We initialize this ensemble by choosing values normally distributed around an estimated value. Here, $\mathbf{x}^{f,(i)}_{t}$ is the state system forecasted from the prior probability distribution (PD) at time $t$, and $\mathbf{P}_{f}$ is the error covariance of that state.\\
\indent Let $\mathbf{y}_{t}$ be an observation at time $t$; the noise on our observations is assumed to be multivariate and normally distributed around $0$. We reflect this using the measurement covariance matrix, $\textbf{R}_t$, i.e. $w^{(i)}_t\sim \mathcal{N}_n(0, \textbf{R}_t)$; let the model noise be multivariate normally distributed around $0$ using the model system covariance matrix, $Q_t$, or $v^{(i)}_t\sim \mathcal{N}_{m_t}(0, \textbf{Q}_t)$; $\textbf{H}_t$ is the observation operator and $\textbf{C}_t=\tilde{\textbf{S}}_t$ where $\tilde{\textbf{S}}_t$ is the sample covariance of the current ensemble. The three equations that describe the analysis step are as follows:

\begin{equation}
\mathbf{x}^{a,(i)}_{t}=\mathbf{x}^{f,(i)}_{t}+ \mathbf{K}^{a}_{t}(\mathbf{y}_{t}-\mathbf{H}^{(i)}_{t}\mathbf{x}^{f,(i)}_{t})
\end{equation}
\begin{equation}
\mathbf{P}^{a}_{f}= \frac{1}{N-1} \sum_{i=0}^N[\mathbf{x}^{f,(i)}_t-\mathbf{x}^a_i][\mathbf{x}^{f,(i)}_t-\mathbf{x}^a_i]'
\end{equation}

\begin{equation}
\mathbf{K}^a_t=\mathbf{C}_t \mathbf{H}'_t(\mathbf{H}_t \mathbf{C}_t \mathbf{H}'_t + \mathbf{R}_t)^{-1}
\end{equation}

\vspace{0.5em}
Here (3) is the state from data posterior PD, where (4) is the covariance matrix for the new prediction, and (5) is the Kalman gain. The notation for the sub-/superscripts are: $f=$ forecast, $a=$ analysis, $t=$ time, and $(i)$ indexes the ensembles. \\
     \indent It is worth noting that often the convention $C_t=Q_t$ is used, and that is used for this research. Further, we calculate the covariance matrix of the analysis with each assimilation, $P = \frac{1}{N-1}\sum_{i=1}^N[x^{(i)}_t-x^a_i][x^{(i)}_t-x^a_i]'$.  $\mathbf{Algorithm}$ $\mathbf{1}$ outlines the algorithm we use for a model in the time range of $t=\{0,1,\dots, T\}$. \\
\begin{algorithm}
    \caption{Ensemble Kalman Filter}\label{alg:alg1}
    \begin{algorithmic}
        \State Generate $\mathbf{x}^{a,(0)}_0, \mathbf{x}^{a,(1)}_0, \dots, \mathbf{x}^{a,(N)}_0$
        \For{$t=0,1,\dots,T$}
            \If{ t in $T_{\text{obs}}$}
                \State Calculate $\mathbf{C}_t$
                \State Calculate $\mathbf{R}_t$
                \State Calculate $\mathbf{K}^a_t=\mathbf{C}_t \mathbf{H}'_t(\mathbf{H}_t \mathbf{C}_t \mathbf{H}'_t + \mathbf{R}_t)^{-1}$
                    \For{ $i=0,1,\dots, N$}
                        \State $\mathbf{x}^{f,(i)}_{t}=\mathcal{M}\mathbf{x}^{(i)}_{t-1}$
                        \State $ \mathbf{y}^{(i)}_t = \mathbf{y}_t+\mathbf{v}^{(i)}_t$
                        \State $\mathbf{x}^{a,(i)}_{t}= \mathbf{x}^{f,(i)}_{t} + \mathbf{K}^a_t (\mathbf{y}^{(i)}_t - \mathbf{H}_t \mathbf{x}^{f,(i)}_{t})$
                    \EndFor
            \State Calculate $\mathbf{x}^a_t = \frac{1}{N} \sum_{i=0}^N \mathbf{x}^{a,(i)}_{t}$
            \State Calculate $\mathbf{P}_f^a = \frac{1}{N-1} \sum_{i=0}^N[\mathbf{x}^{f,(i)}_t-\mathbf{x}^a_i][\mathbf{x}^{f,(i)}_t-\mathbf{x}^a_i]'$
            \Else{}
                \For{ $i=0,1,\dots, N$}
                    \State $\mathbf{x}^{(i)}_{t}=\mathcal{M}\mathbf{x}^{(i)}_{t-1}$
                \EndFor
                \State Calculate $\mathbf{x}^a_t=\frac{1}{N}\sum_{i=0}^N\mathbf{x}^{a,(i)}_{t}$
                \State Calculate $\mathbf{P}_f = \frac{1}{N-1}\sum_{i=0}^N[\mathbf{x}^{(i)}_t-\mathbf{x}^a_i][\mathbf{x}^{(i)}_t-\mathbf{x}^a_i]'$
            \EndIf
    \EndFor
    \end{algorithmic}
\end{algorithm}
\indent Can data assimilation also provide more accurate and timely forecasts of glacier melt? What is the computational cost, and can our experiments with data assimilation determine in what time frame observational data would be most critical to the quality of the model output?

\subsubsection{Twin Experiments}
\indent We work with observing system simulation experiments (OSSE), or twin experiments, to perform different ``stress tests" of our data assimilation system within a well-structured framework that does not require us to work with real-world data. The method is as follows: We fix all parameters and unknowns and define from the direct model a run which we call the ``truth". We perturb a set of synthetic measurements using a random normal distribution from the ``true" run to generate our ``observations". Then we perturb our initial conditions and unknowns and define a false ``initial guess", different from the ``true" initial state, which we are then able to use to evaluate the performance of the method. This is a standard tool for evaluating inverse problems (\cite{asch2016fundamentals}).

\subsubsection{Square Difference}\label{sec:sqdf}
\indent The error measure we use to determine the best ensemble size and observation scheme (Figure \ref{fig:ensvL}) is the square difference, $d^2$, which we define as
\begin{equation} d_t^2 = \left(x_t - x^a_t\right)^2\end{equation}
where $x_t$ is the true state from the truth simulation at time $t$ and $x^a_t$ is the analysis state at time $t$.

\subsection{Sensitivity Analysis}
\indent We also conduct sensitivity analyses to understand how various sources of uncertainty in a mathematical model contribute to the model's overall uncertainty (\cite{modelanalysis}). In our analysis, we vary parameters by +/-10 percent of the nominal values originally given in our model code, with the understanding that if the outputs vary significantly, the output is considered sensitive to the specification of the input distributions. This provides insight into the most promising applications of data assimilation. For this analysis, the parameters are grouped together into three categories: initial conditions (initial height, length, and bed slope), sill parameters (sill minimum, maximum, and slope), and smb parameters (initial, mid, and final smb). The parameters in each category are perturbed together for simplicity, as they are related to one another. We analyze the spread of the glacier heights and lengths 
produced from the range of perturbed values for each of the three categories. 

\section{Results}
In the following sections, we present the effect the data assimilation scheme has on the forecasts of glacier melt.

\subsection{Sensitivity Analysis}
\begin{figure*}[htbp]
   \subfloat[$\pm 10 \%$ surface mass balance variables]{\label{fig: t vs H(t) varying SMB}
  \label{fig:smb}
      \includegraphics[width=.32\textwidth]{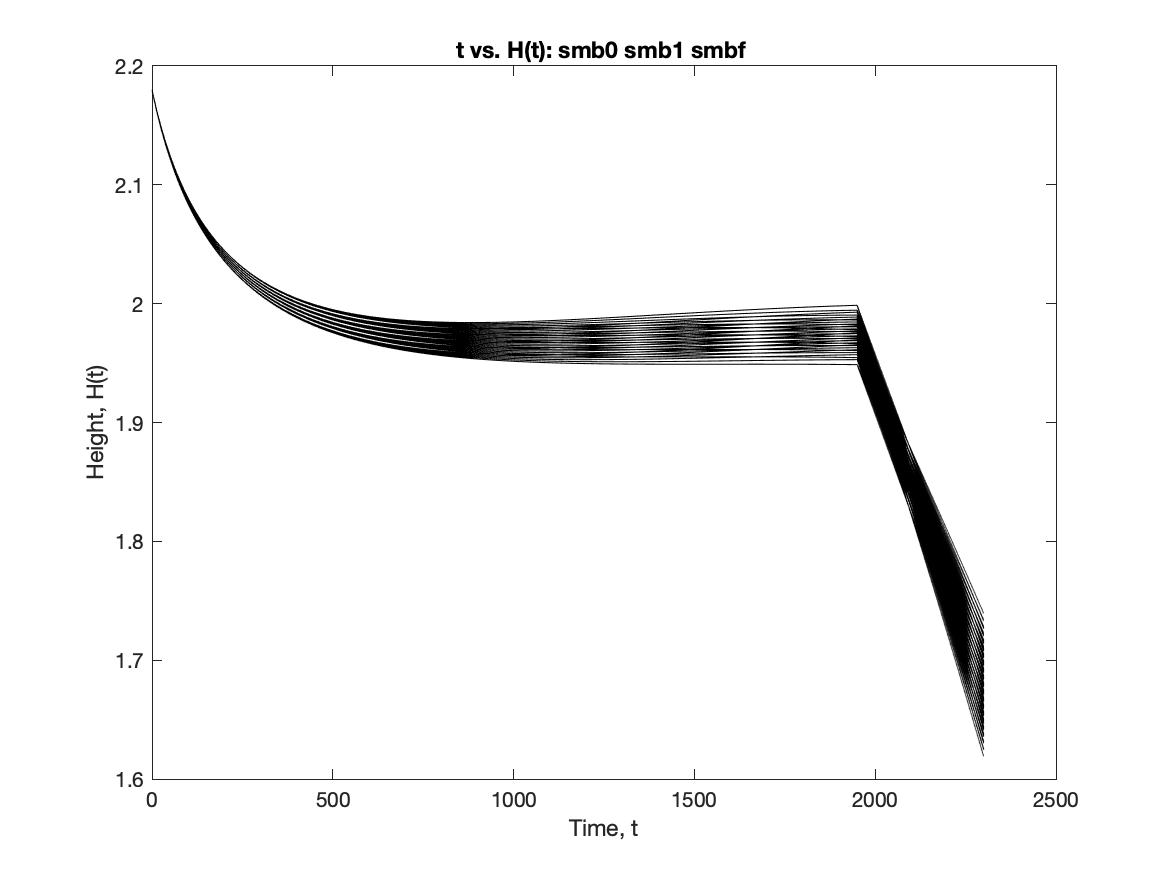}}
~
   \subfloat[$\pm 10 \%$ initial conditions and bed slope]{\label{fig:t vs H(t) varying initial conditions}
  \label{fig:init}
      \includegraphics[width=.32\textwidth]{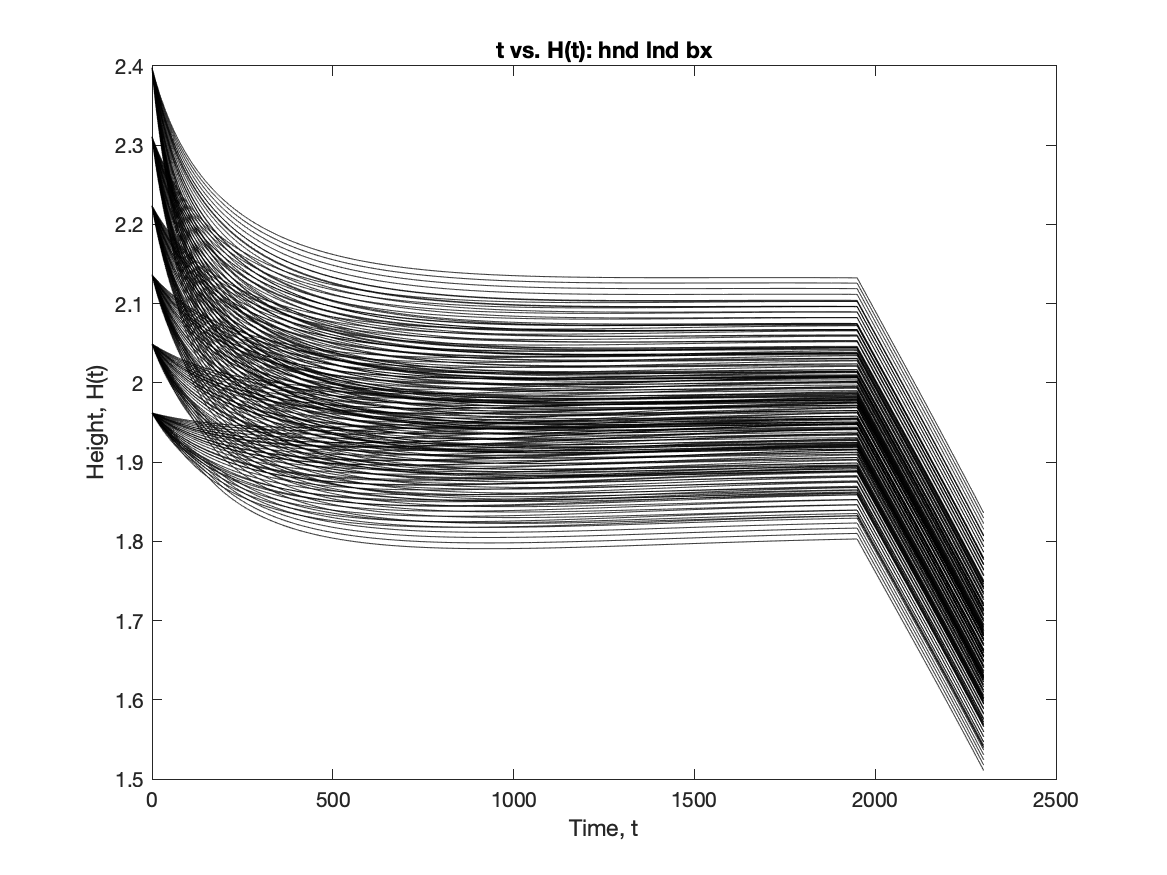}}
~
   \subfloat[$\pm 10 \%$ sill variables]{\label{fig: t vs H(t) varying sill}
  \label{fig:sill}
      \includegraphics[width=.32\textwidth]{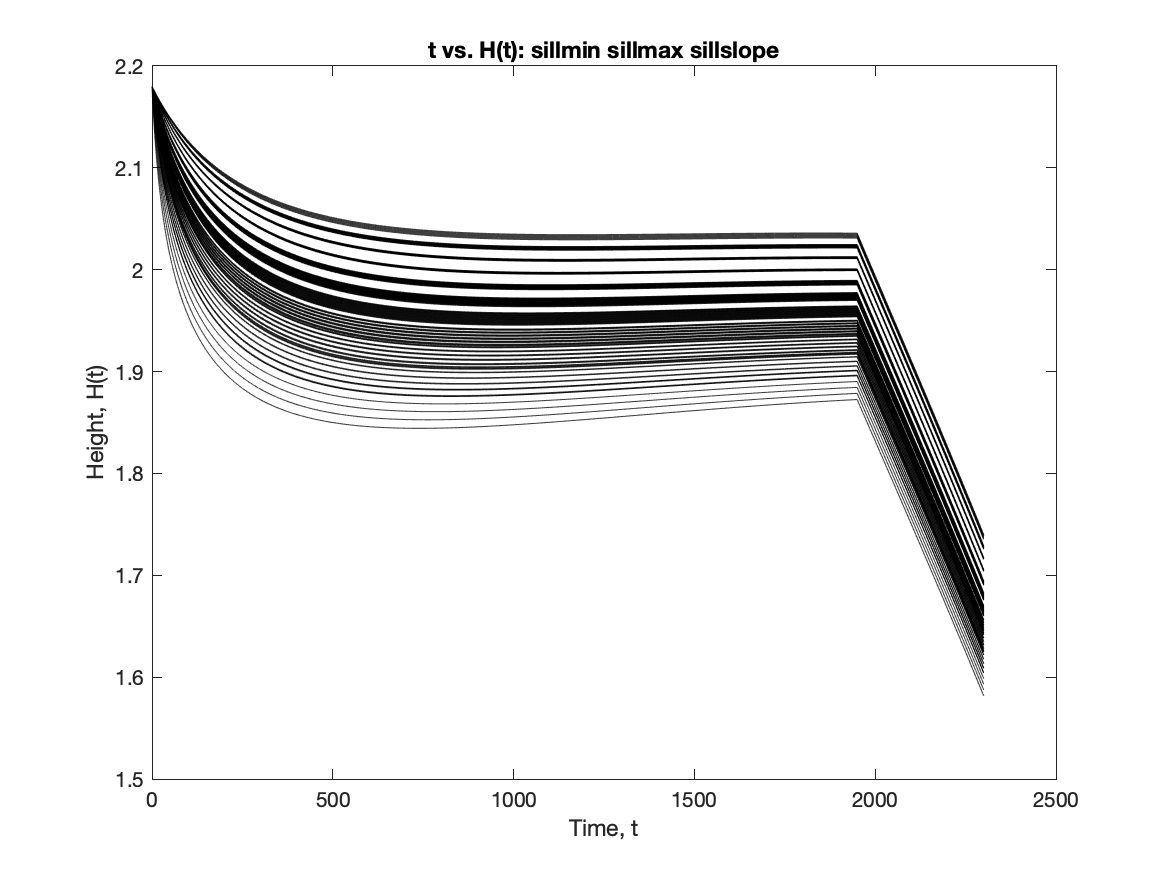}}

   \caption{Plots of $t$ vs $H(t)$ varying SMB, initial conditions, and sill values respectively}\label{bs1}
\end{figure*}

\indent In our study of the effect of varying the parameters of the model by $\pm$ 10\%, we see that in Figures \ref{fig:smb}, \ref{fig:init} and \ref{fig:sill} that the height of the glacier decreases until 1950, then decreases more steeply until 2300 for all parameter values simulated. Despite this general pattern, we see that this implementation of 10\% variation causes both a great deal of variation in model behavior and very little.  

\begin{wrapfigure}{r}{0.5\textwidth}
    \centering
    \includegraphics[width=0.9\linewidth]{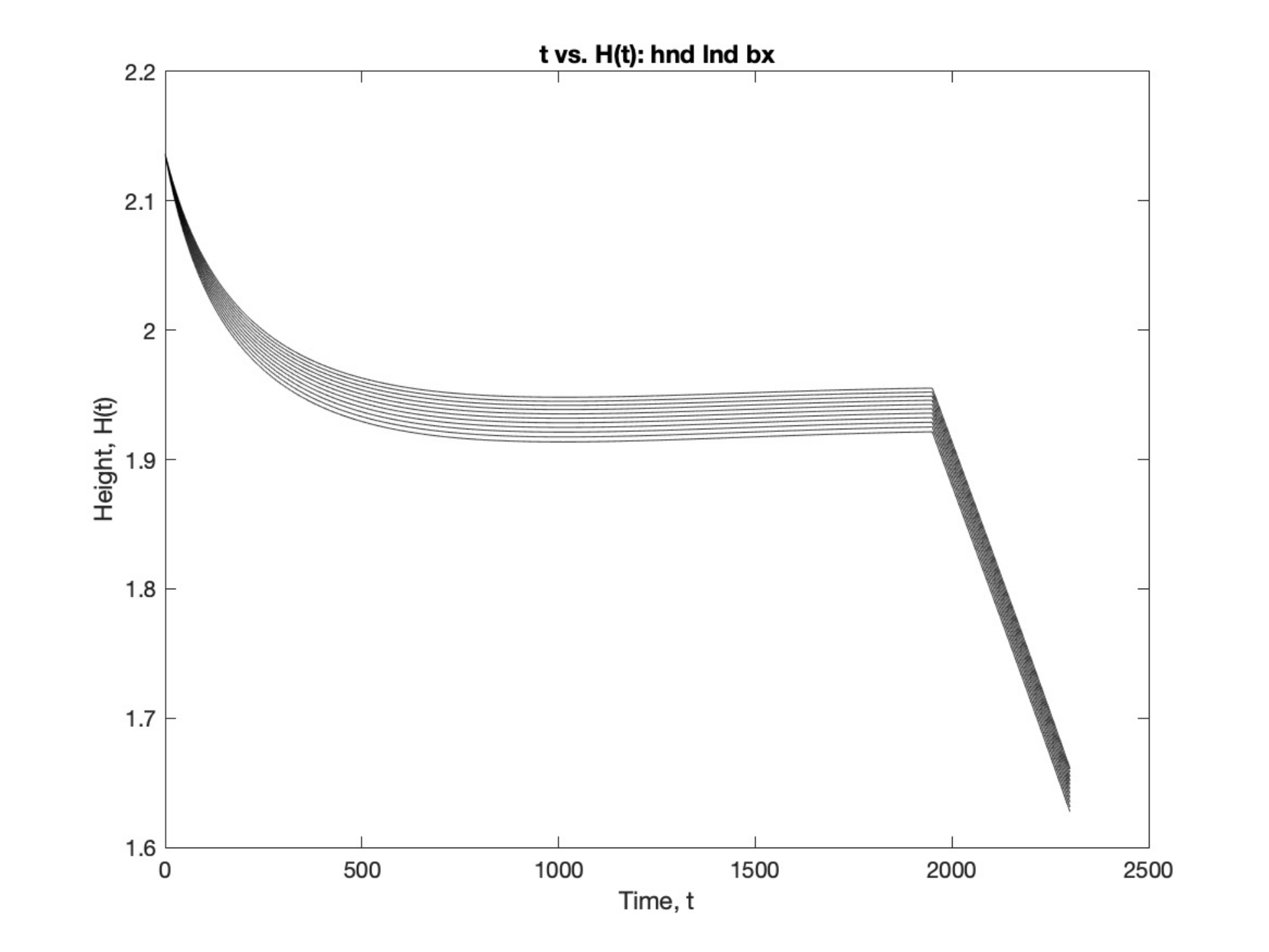}
    \caption{Single initial condition with slope varied $\pm 10\%$.}
    \label{fig:singlecase}
\end{wrapfigure}

\indent As we can see in Figure \ref{fig:smb}, varying the SMB does not cause very much variation in the  model run, and all simulations produces a generally similar shape.
Varying initial conditions along with slope \ref{fig:init} and the sill variables \ref{fig:sill} causes a great deal of variation. However, it should be noted that much of the variation in Figure \ref{fig:init} is actually due to changing the initial condition we chose in Figure \ref{fig:initcondtrue}, as we will see if we vary slope for a single initial condition as in Figure \ref{fig:singlecase}.
Finally, we see that for the sill variables in Figure \ref{fig:sill} that changing their values caused a great deal of variation. While each run maintained the same general shape, we see that the curvature before 1950 

varies from a strong dip to seemingly approaching a constant before 1950.\\
\indent Thus, when we decide on conditions and parameters to model a "truth" simulation, we know that if we change the sill parameters we get a more deviated case (and simulate real world measurement errors) to assimilate back to our truth simulation. We discuss our choice of initial conditions in the next section.

\subsubsection{Parameters and Initial Conditions for Our Study}\label{ourParameters}
\indent Presented in the following table are the parameters and initial conditions used for our study. The ``True" column holds the values for the truth simulation, and the ``Inaccurate" column holds the slightly perturbed values that simulate our initial assumptions of truth values which would be hidden from us in the study of a real glacier. 
\begin{figure}[H]
\centering
\begin{tabular}{|p{1cm}||p{1cm}|p{1.6cm}|p{2.3cm}|  }
 \hline
 \multicolumn{4}{|c|}{Parameters and Initial Conditions} \\
 \hline
 & True & Inaccurate & \% variation\\
 \hline
$smb_o$ & 0.3 & 0.35 & $+\approx$ 16.67 \\
$smb_1$ & 0.15 & 0.15 & n/a\\
$smb_f$ & 0.0 & 0.0 & n/a \\
$H_o$ & 2.18 & 2.3 & $+\approx$ 5.5\\
$L_o$  & 4.44 & 4.6 & $+\approx$ 3.3\\
$b_x$  & -0.001 & -0.001 & n/a\\
$sill_{min}$ & 415 & 415 & n/a\\
$sill_{max}$ & 425 & 425 & n/a\\
$sill_{slope}$ & 0.01 & 0.008 & -20\\
 \hline
\end{tabular}
\end{figure}

\subsection{Data Assimilation}

\subsubsection{Ensemble Size}\label{ensembleSize}
\begin{figure*}[htp]
   \subfloat[]{\label{fig:ensvH}
      \includegraphics[width=.46\textwidth]{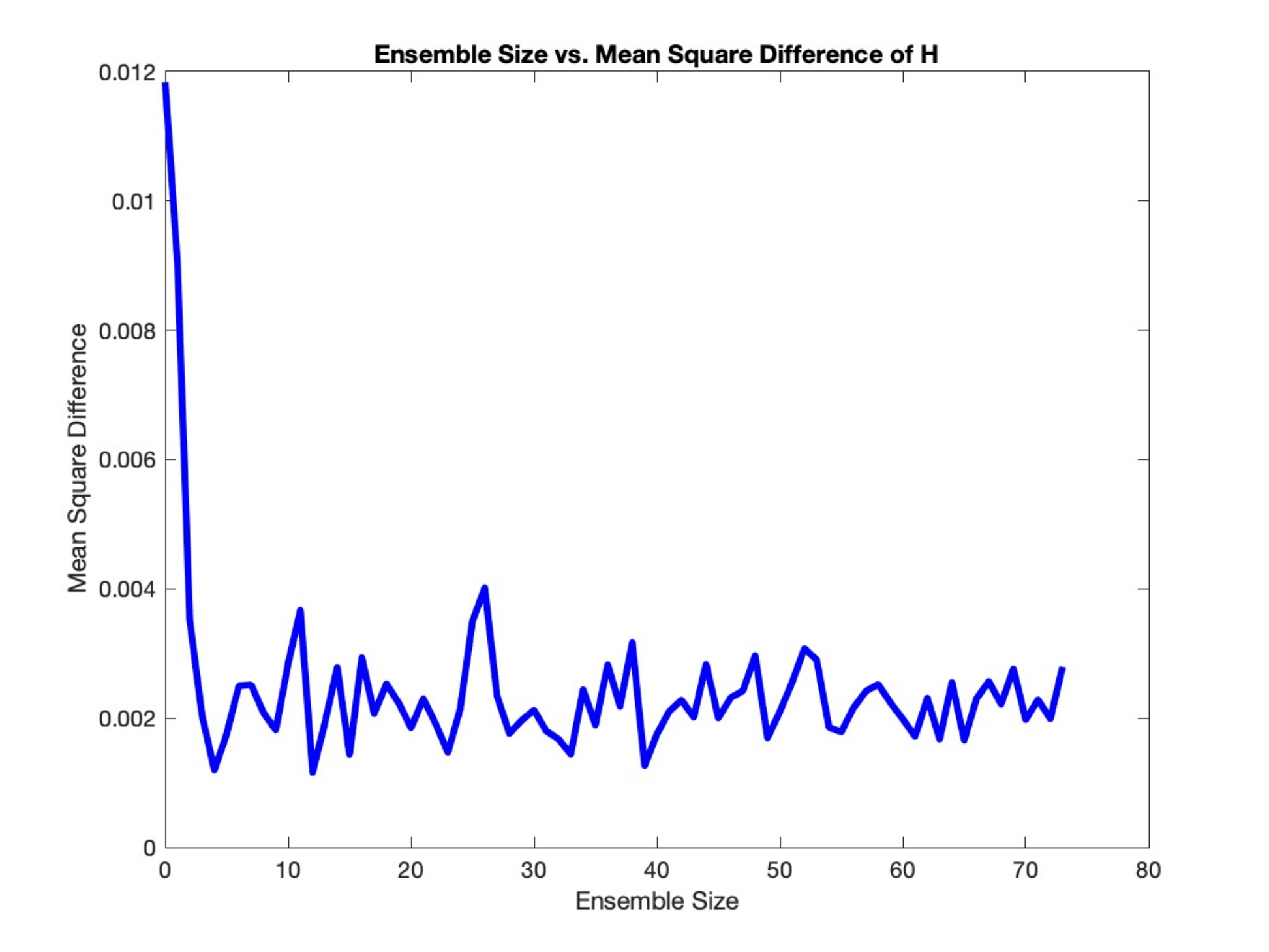}}
~
   \subfloat[]{\label{fig:ensvL}
      \includegraphics[width=.46\textwidth]{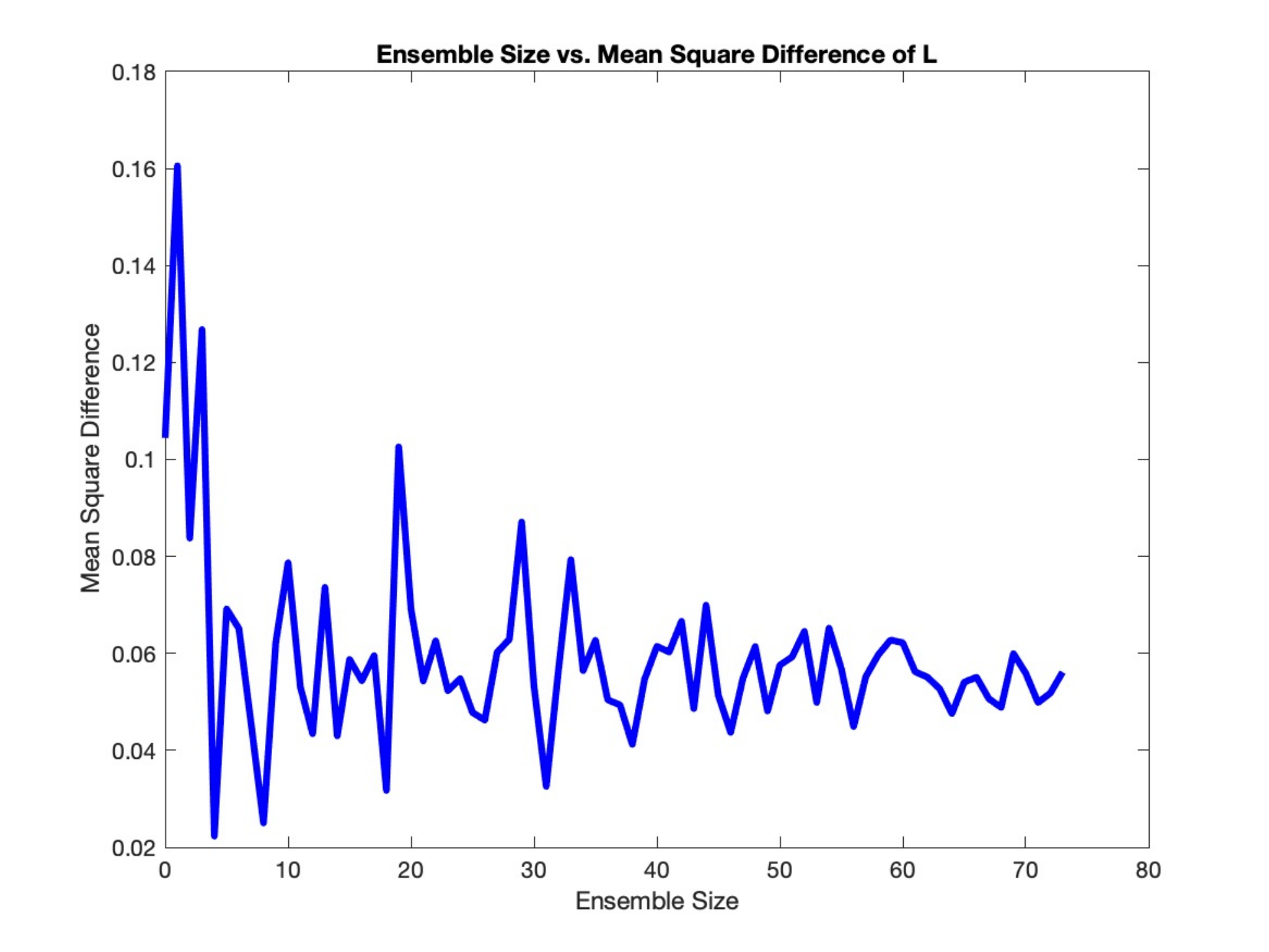}}

   \caption{Mean Square Difference for each Ensemble Size in height H (a) and length L (b)}
\end{figure*}

\indent Using the square difference as described in Section \ref{sec:sqdf}, we examine the mean square difference for various ensemble sizes. We ran this calculation for ensembles sizes ranging from 2 to 75, and found that ensembles of sizes 7-10 were ideal, as the mean square difference hovers around the same value with larger ensemble sizes (Figures \ref{fig:ensvL}).

\subsubsection{Observation Scheme}\label{observationScheme}
\indent We run the model for various observation schemes to find the optimal scheme, i.e. the times frames and frequencies which can produce a sufficiently small average square difference over the course of the model run. We applied this process to our model and found that for before 1900 the best observation frequency, while still using small number of observations, would be every 19 years for a total of 100 observations (Figure 9). Similarly, for the time frame of 1950-2300 we found that yearly observations for a total of 350 observations is the best frequency (Figure 10). Refer to Appendix ~\ref{fig: assimilationDates} for a graphical comparison of assimilation frequency runs.

\begin{figure*}[htp]
   \subfloat[Assimilated every 19 years]{\label{fig:oldData19Years}
      \includegraphics[width=.46\textwidth]{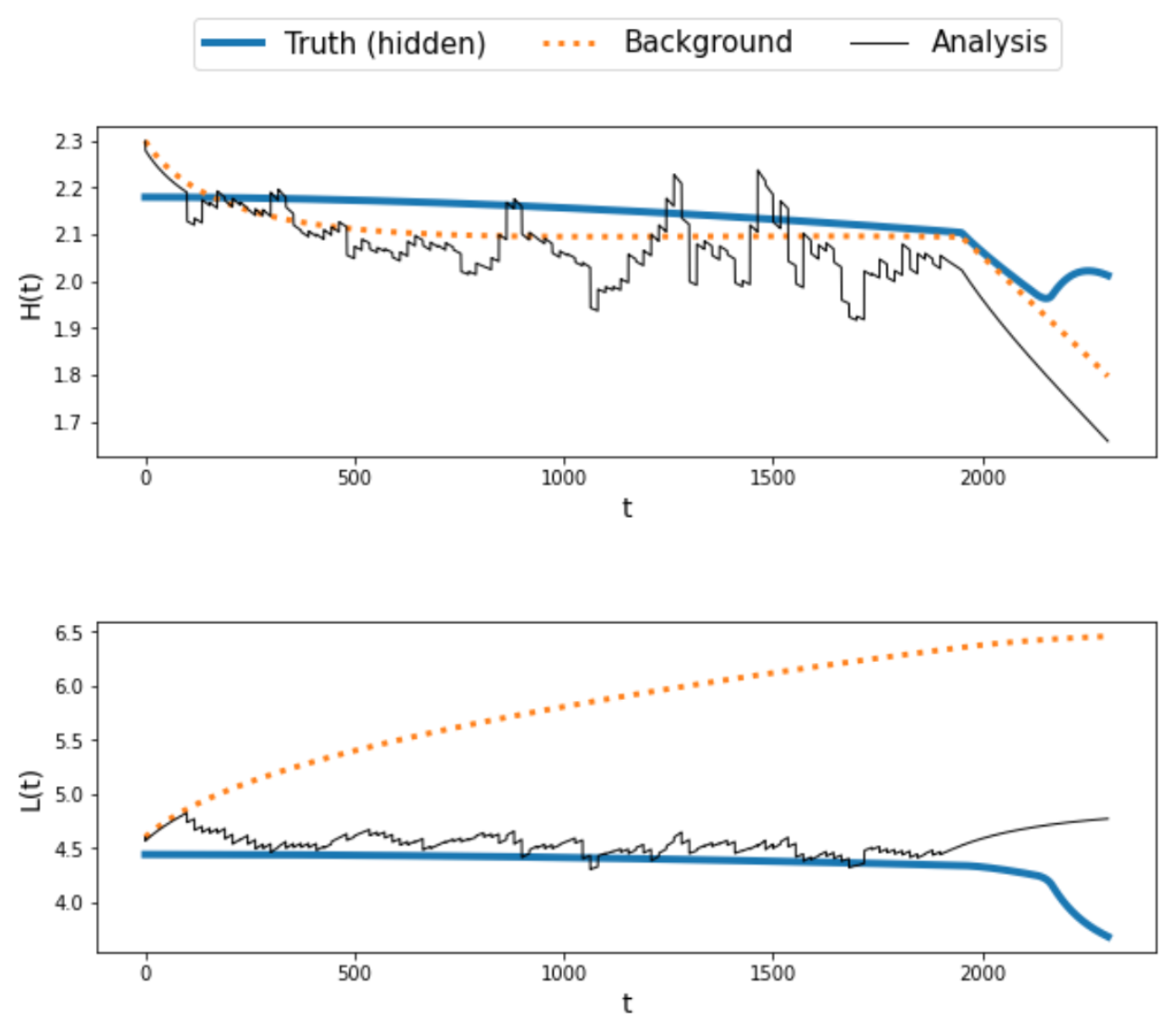}}
~
   \subfloat[Assimilated yearly]{\label{fig:newDataYearly}
      \includegraphics[width=.46\textwidth]{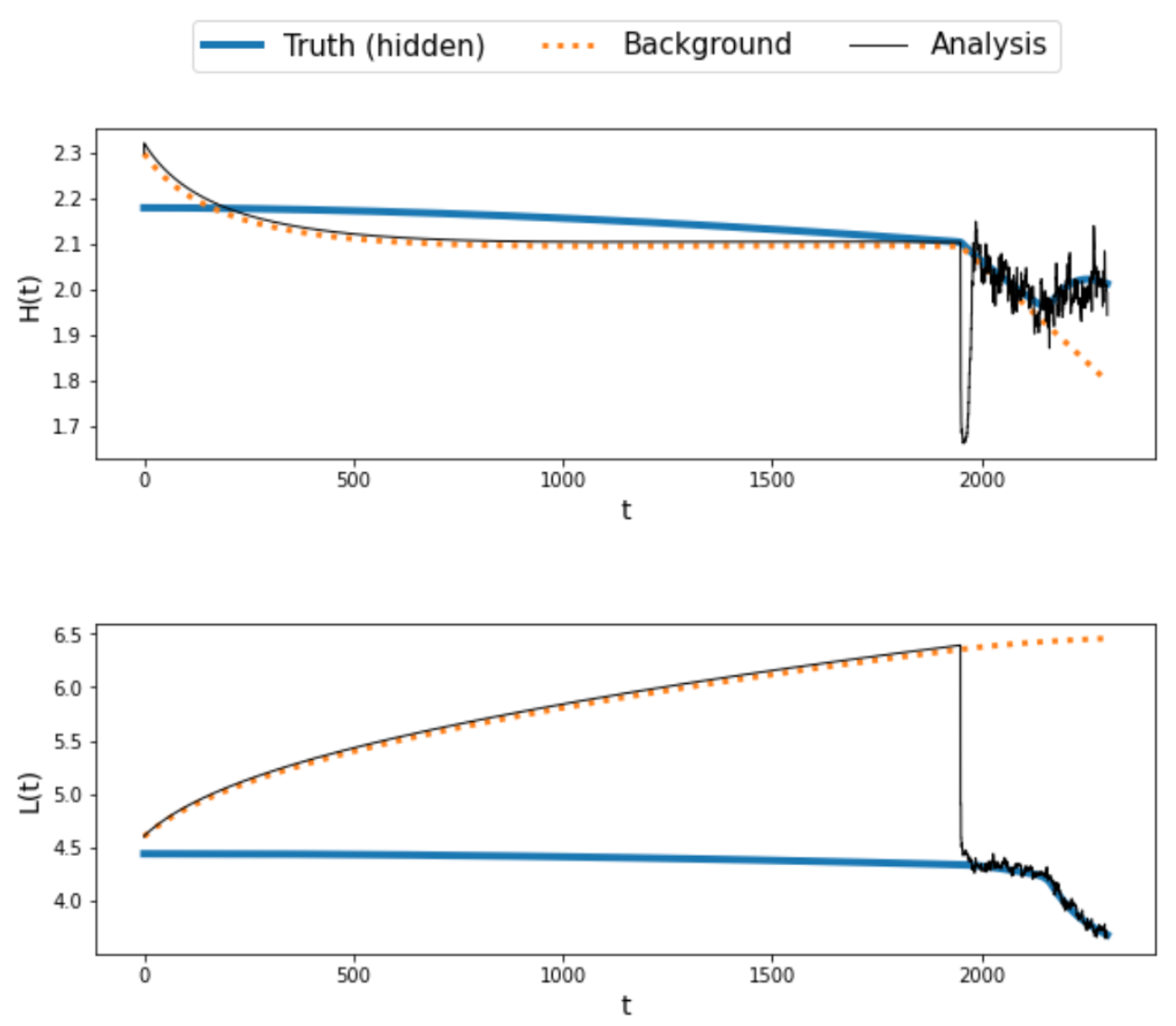}}

   \caption{Background, truth, and analysis for old data assimilated every 19 years (a) and assimilated yearly (b)}
\end{figure*}

\begin{figure*}[htp]
   \subfloat[Assimilated every 19 years]{\label{fig:oldDataMSE}
      \includegraphics[width=.46\textwidth]{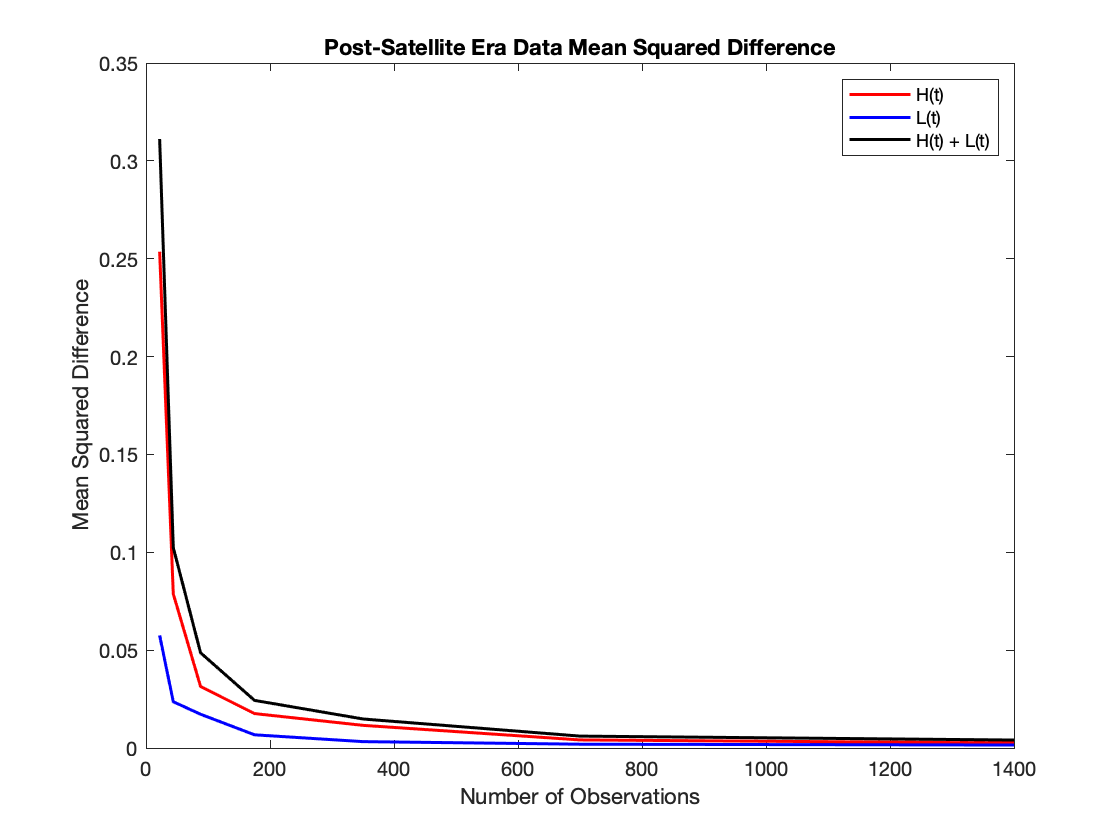}}
~
   \subfloat[Assimilated yearly]{\label{fig:newDataMSE}
      \includegraphics[width=.46\textwidth]{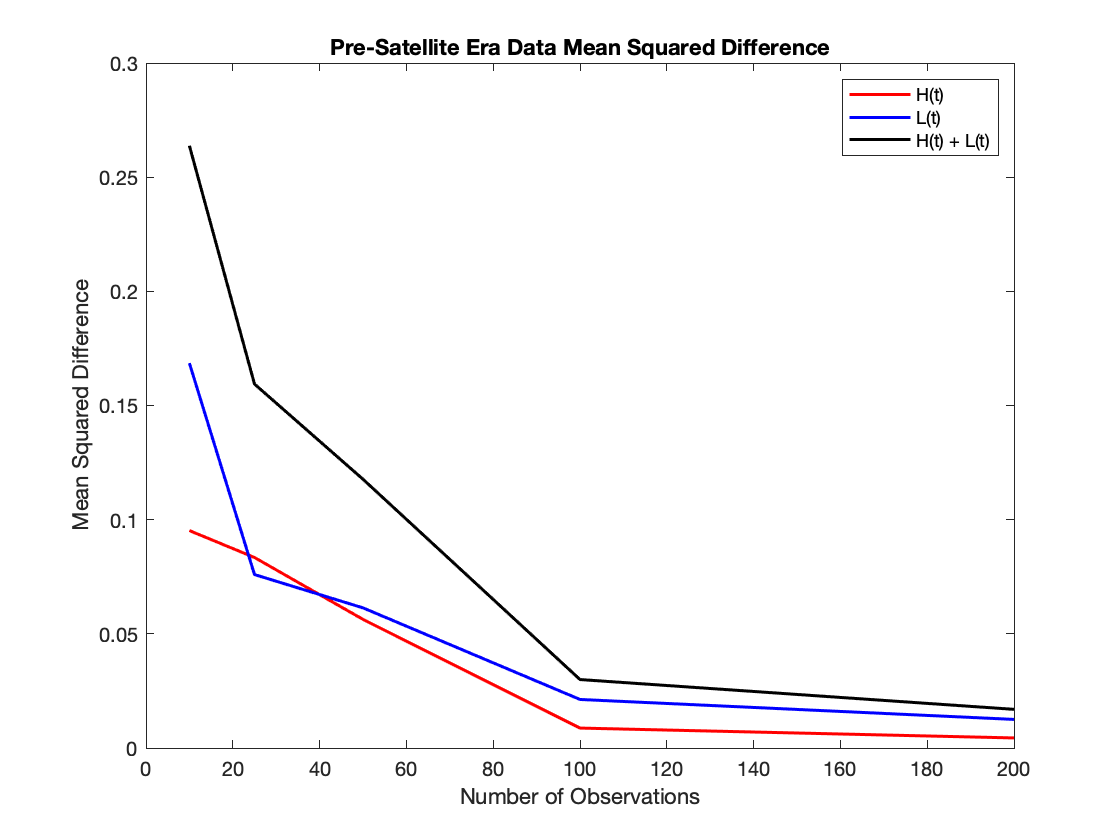}}

   \caption{Mean squared difference between analysis and truth for pre-satellite era data (a) and post-satellite era (b) data assimilation with different frequencies}
\end{figure*}

\subsection{Model Runs}
\indent We find that the best set up for the EnKF is an ensemble size of 7-10, pre-1900 observation frequency of every 19 years and yearly observations after 1950. Using these facts, we ran the model using EnKF for the time frame of 0-2300 using the ``best" observation and ensemble scheme. For comparison, we also ran the model using EnKF for the time frame of 0-2300 using a ``worse" observation and ensemble scheme with an ensemble size of 10, observations every 200 years starting in year 200, and measurements every year starting in 2000. We found that  both runs were generally accurate, but the best observation and ensemble scheme was an improvement over our worse case. \\
\indent Next, we ran the model following the best observation and ensemble scheme except with observations terminating in 2022 in order to project $H$ and $L$ into the future up to the year 2300. The following plots show the results of this experiment, which we will use to help calculate $Q_g$ over time, and in turn use it to calculate sea level rise.
\begin{figure*}[htp]
   \subfloat[]{
      \includegraphics[width=.46\textwidth]{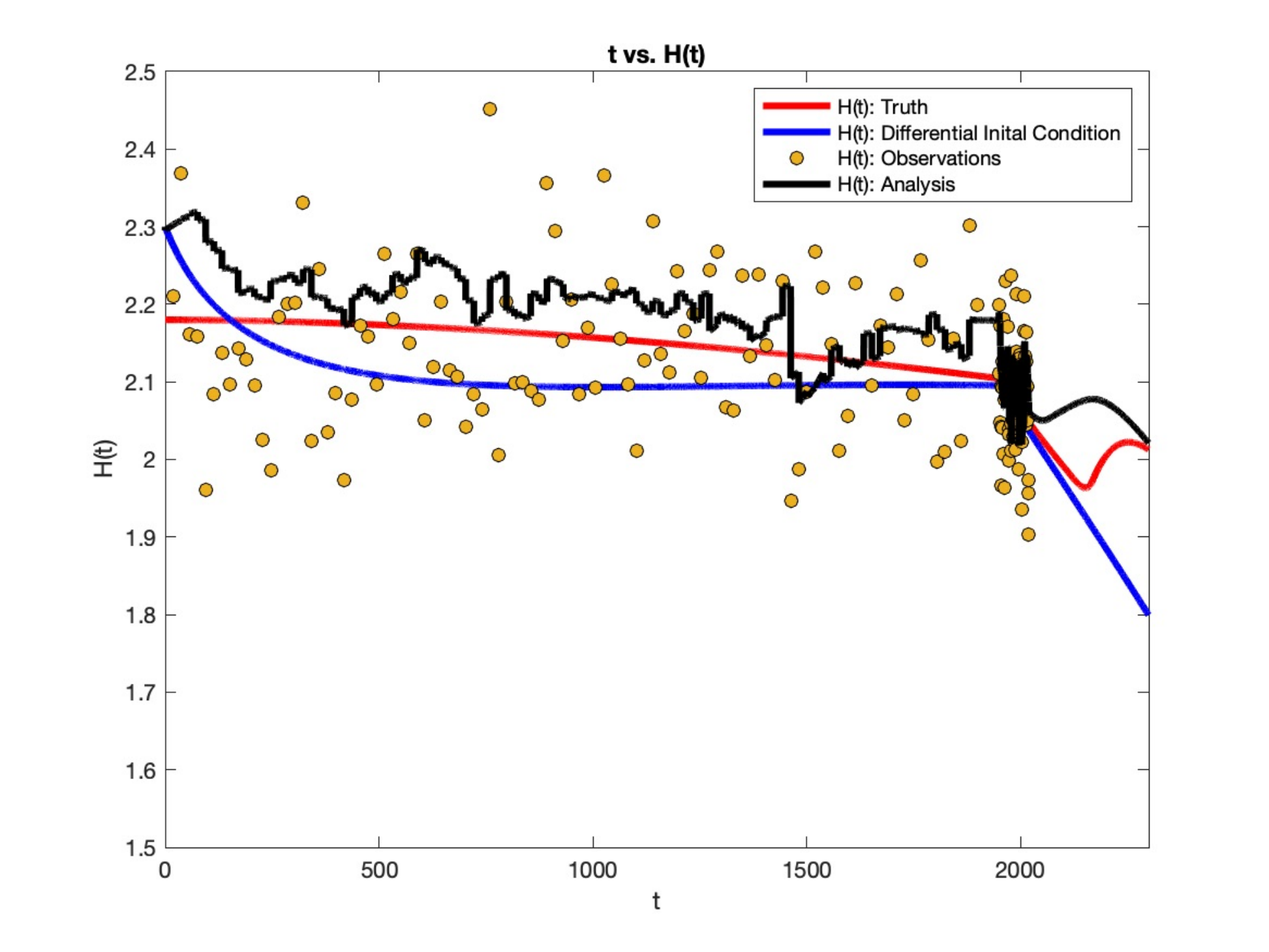}}
~
   \subfloat[]{
      \includegraphics[width=.46\textwidth]{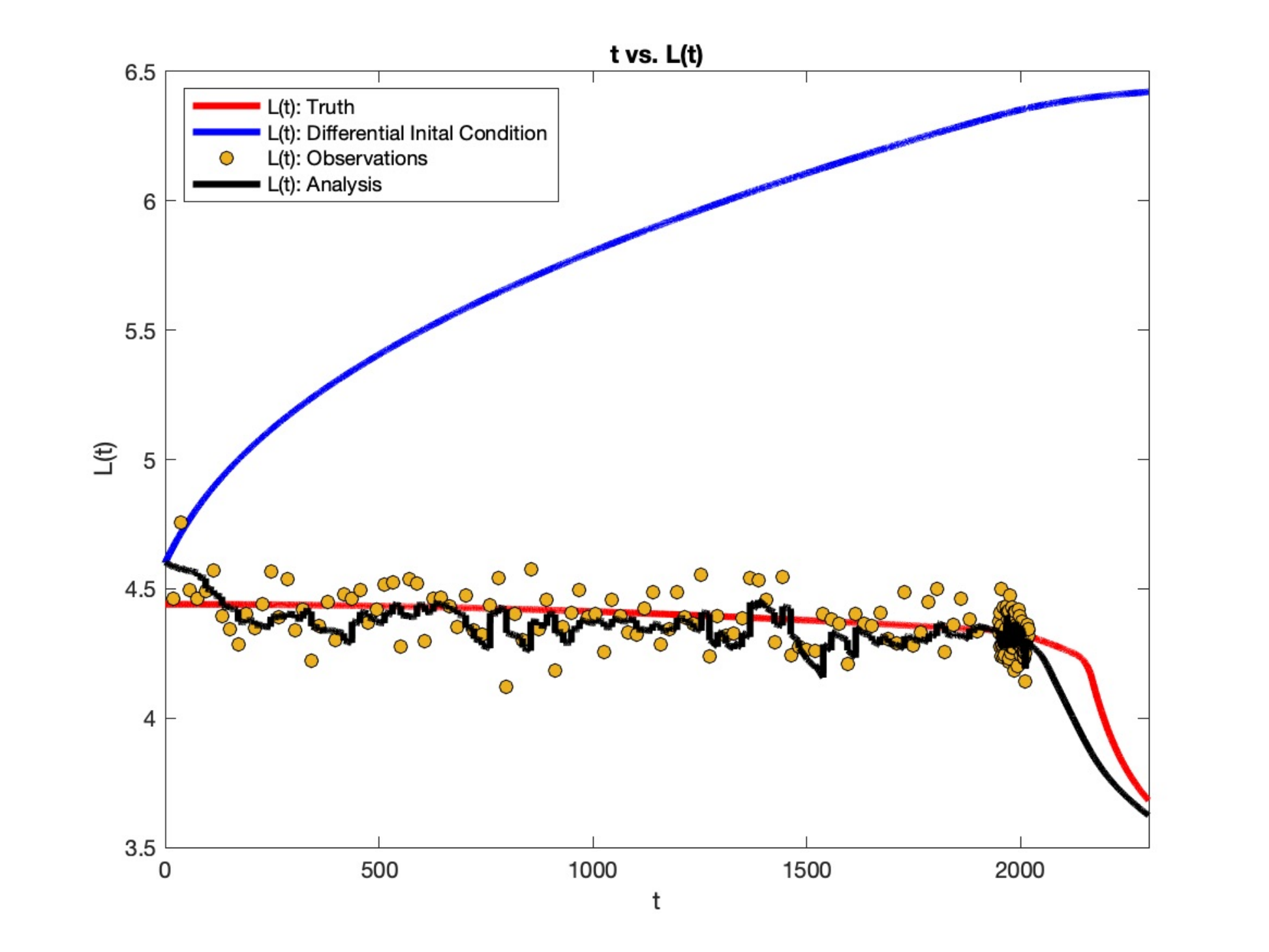}}

   \caption{Model run of H (a) and L (b) using best observation scheme and ensemble size. Note that we only have our final analysis displayed in the plot}
\end{figure*}

\begin{figure}[H]
\centering
\begin{tabular}{|p{1.5cm}||p{1.5cm}|p{1.5cm}|  }
 \hline
 \multicolumn{3}{|c|}{Projected Values} \\
 \hline
Year & H & L \\
 \hline
2000  & 2.1672 & 4.3532 \\
2050 & 2.028 & 4.3028 \\
2100 & 2.0207 & 4.1973 \\
2150 & 2.0442 & 4.0063 \\
2200 & 2.0549 & 3.845 \\
2250  & 2.0461 & 3.7271 \\
2300  & 2.0219 & 3.6438 \\

 \hline
\end{tabular}
\caption{Values of the model at various times. Note that the year 2000 was not a projection from our model but just put as a reference.}
\end{figure}

\subsection{Sea Level Rise} 
As this glacier model calculates the behavior of the average length, $L$, and average thickness, $H$, over time, we can see how this process would prompt questions about the glacier mass lost after crossing the grounding line into the sea. 

\begin{figure}[H]
        \centering
        \scalebox{0.8}{\begin{tikzpicture}[
        squarednode/.style={rectangle, draw=blue!60, fill=blue!5, very thick,
        minimum size=5mm}]
        \node[squarednode]      (middle)                     {Sea Level Rise};
        \node[squarednode]      (leftsquare)       [left=of middle] {Height and Length of Glaciers};
        \node[squarednode]        (rightsquare)       [right= of middle] {Storm Surge};
        
        \draw[-latex] (leftsquare.east) -- (middle.west);
        \draw[-latex] (middle.east) -- (rightsquare.west);
        \end{tikzpicture}}
        \caption{Relationship of glacier melt to storm surge}
        \label{fig:seaLevelImpactDiagram}
    \end{figure}
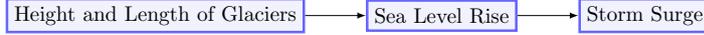

Recall the grounding line and interior fluxes, $Q$ and $Q_g$, defined in Equations \ref{q} and \ref{qg}.
We can use Equations \ref{q} and \ref{qg}, and the following equation from Robet et. al 2019 to calculate the volume lost across the grounding:
\begin{equation}
\dfrac{dV_{gz}}{dt} = W(Q-Q_g).
\end{equation}
\indent If we integrate this equation assuming its initial condition is 0, we find
\begin{equation}V_{gz} = W(Q-Q_g)t.\end{equation}
We use this calculation to estimate the volume lost at each time step and summed them 
\begin{wrapfigure}{l}{0.45\textwidth}
    \centering
    \includegraphics[width=\linewidth]{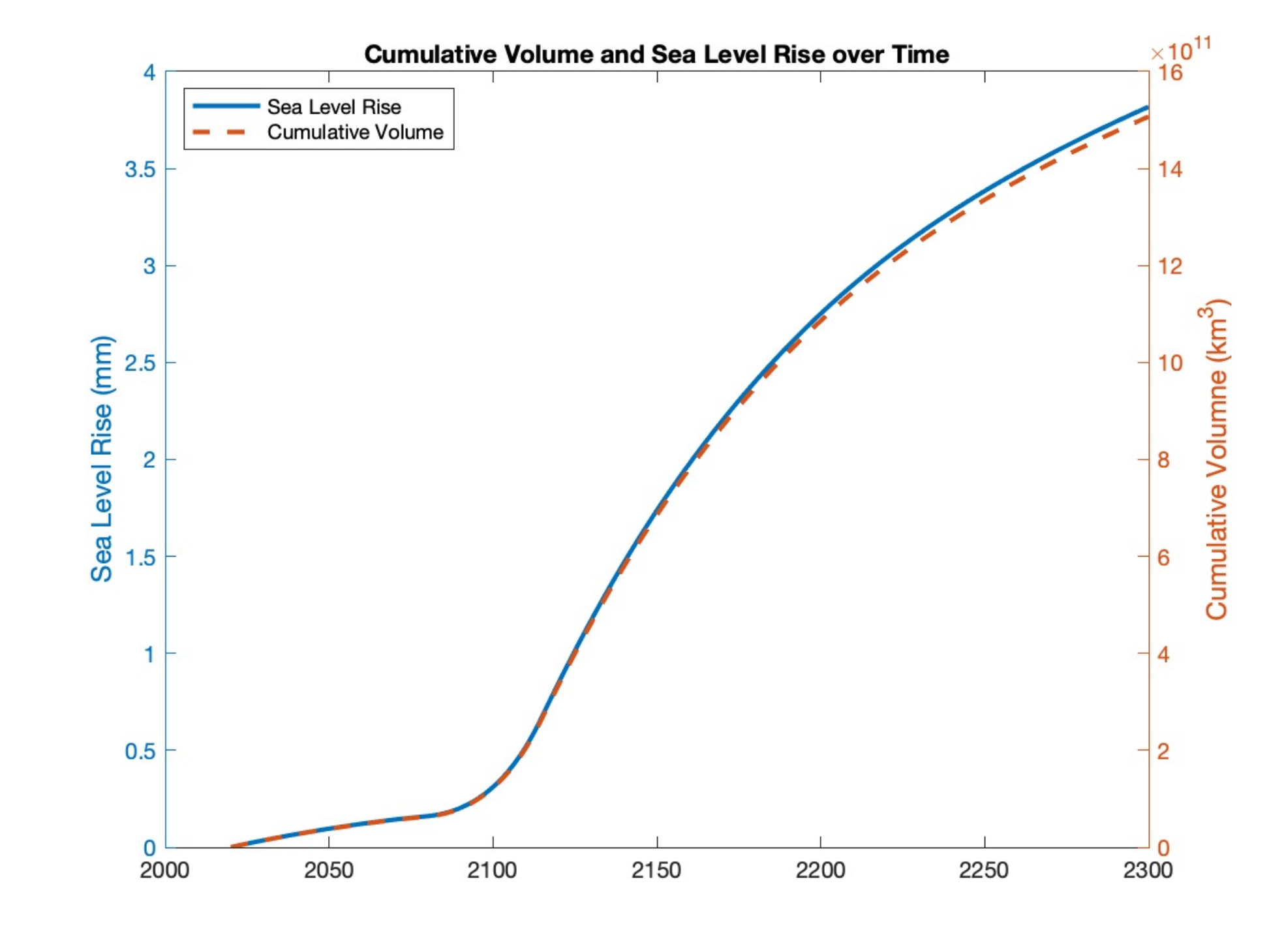}
    \caption{Sea Level Rise and Cumulative Volume for width of 50 km}
    \label{fig:volumesealevel}
\end{wrapfigure}
to estimate total accumulated volume loss. Because we cannot calculate $W$ from our model, we perform this calculation for glaciers of width $5$ km, $50$ km, and $100$ km.
To convert this to sea level rise, note that 394.67 km$^3$ of ice is equivalent $1$ mm of sea level. Notice that sea level rise and cumulative volume out of glacier as shown in Figure \ref{fig:volumesealevel} has the same shape but different values, which is to be expected since we divide by the conversion value to get sea level rise. 
We see that for widths between 5 km and 100 km, we get sea level rise ranging from 0.345 mm to 6.896 mm in the year 2300. We can add these values to the current sea level in order to project the overall future sea level. 


In order to get an approximation for sea level rise we use the fact that there are 733 glaciers in Greenland according to \cite{bjork2015brief}. In order to get an approximate value for total sea level rise, we will crudely assume the width of all the glaciers are the same and multiply its sea level rise by 733.

\subsection{Storm Surge}
\begin{wrapfigure}{r}{0.6\textwidth}
    \includegraphics[width=0.9\linewidth]
{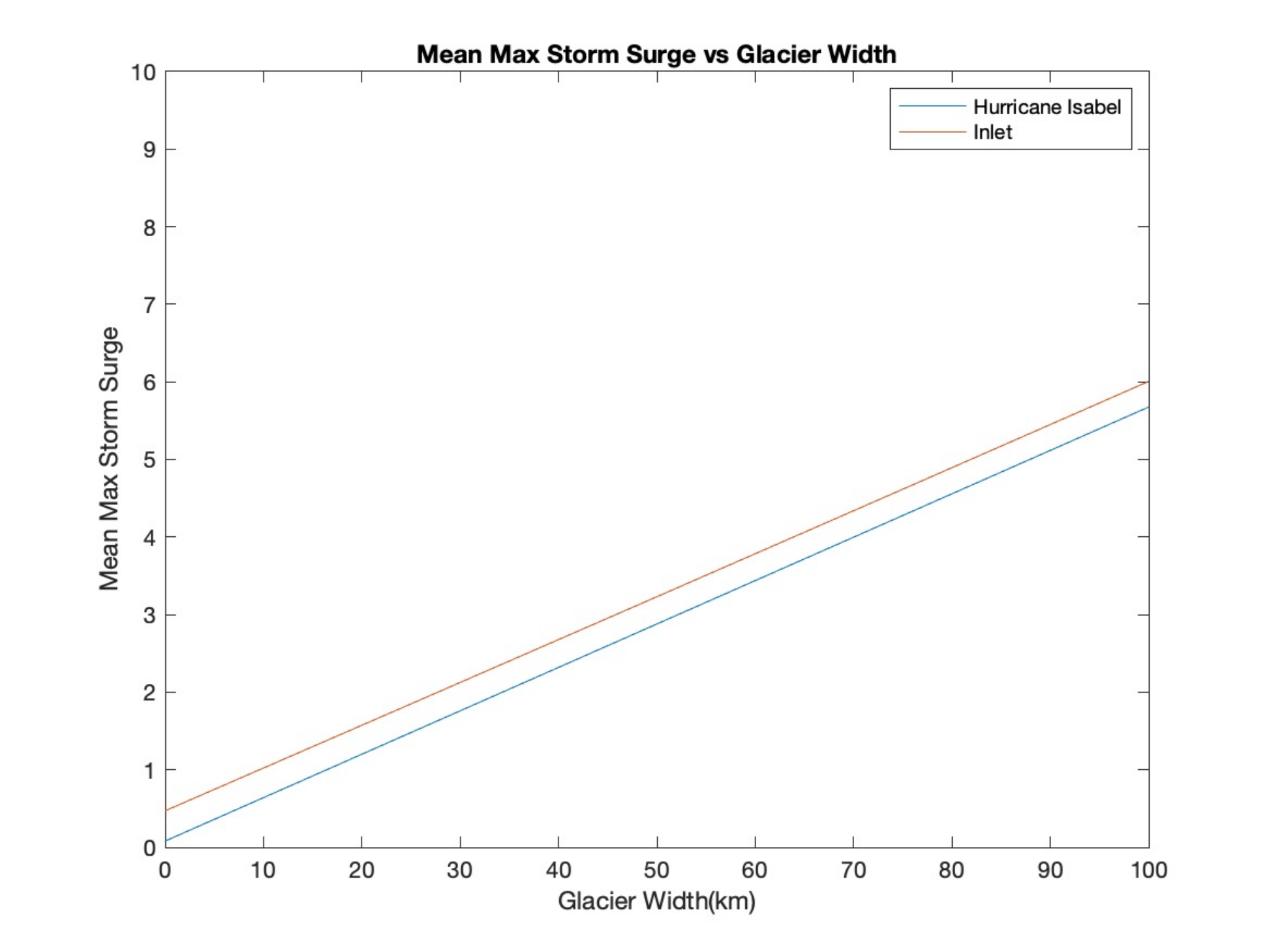}
    \caption{Mean Maximum Storm Surge at a location for the respective simulated Glacier Width}
    \label{fig:stormsurge}
\end{wrapfigure}
We model storm surges using the Advanced Circulation (ADCIRC) model assume these different values of total sea level rise. We run this model for various glacier widths, and determine the sea level rise that these are projected to cause in our model by the year 2300. We then look at differences in the max storm surge and the average change in storm surge over time. To do this, we look at the average max storm surge elevation across all grid points for two example problems from the ADCIRC website: Idealized Inlet and Hurricane Isabel. We see across all points in the model which had storm surge that the mean of the maximum storm surges increases when modeled for greater glacier width alongside the subsequent sea level rise by 2300. For the 5 km cases, we see that the Idealized Inlet increased by 445.6\% and Hurricane Isabel increased by 158.6\%. Thus, we find that there is a direct, positive relationship between the sea level rise caused by glacier melt and storm surge.

\begin{wrapfigure}{l}{0.6\textwidth}
    \includegraphics[width=0.9\linewidth]
{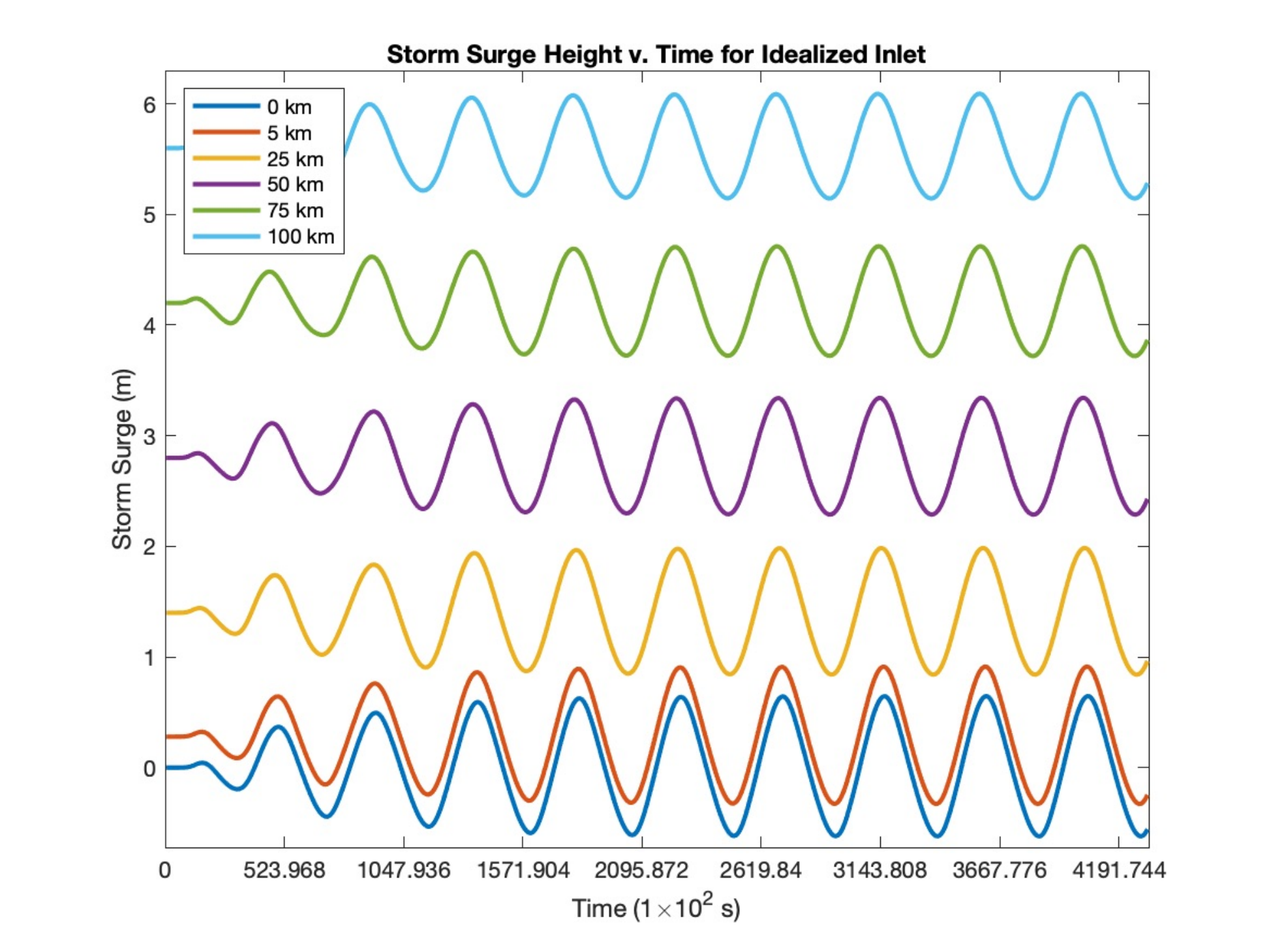}
    \caption{Storm Surge Height over Time for the Idealized Inlet given various Glacier Widths}
    \label{fig:inlettime}
\end{wrapfigure}

\indent Further, when we choose a mesh grid point in the idealized inlet simulation, we see that storm surge increases directly with greater sea level rise(Figure \ref{fig:inlettime}). On the storm surge for the 6 different sea level rise situations, we see that the behavior of storm surge level is similar. Specifically, we observe a slight decrease in magnitude of change for higher sea level rise.

\indent Now that we have seen that there appears to be an overall rise in storm surge levels with sea level rise, we want to check if this is still the case for a simulation of a real life event. To look at the affect of sea level rise on storm surge levels, we look at the square difference between the maximum storm surge over time at each grid point for the 0 km width and the 25 km width simulations. We see that the largest square difference is in the Chesapeake Bay Region, with the square difference being similar everywhere else. The similar square difference is likely a product of general sea level rise (it is constant throughout the model) and not storm surge.

\begin{wrapfigure}{r}{0.6\textwidth}
    \includegraphics[width=0.9\linewidth]
{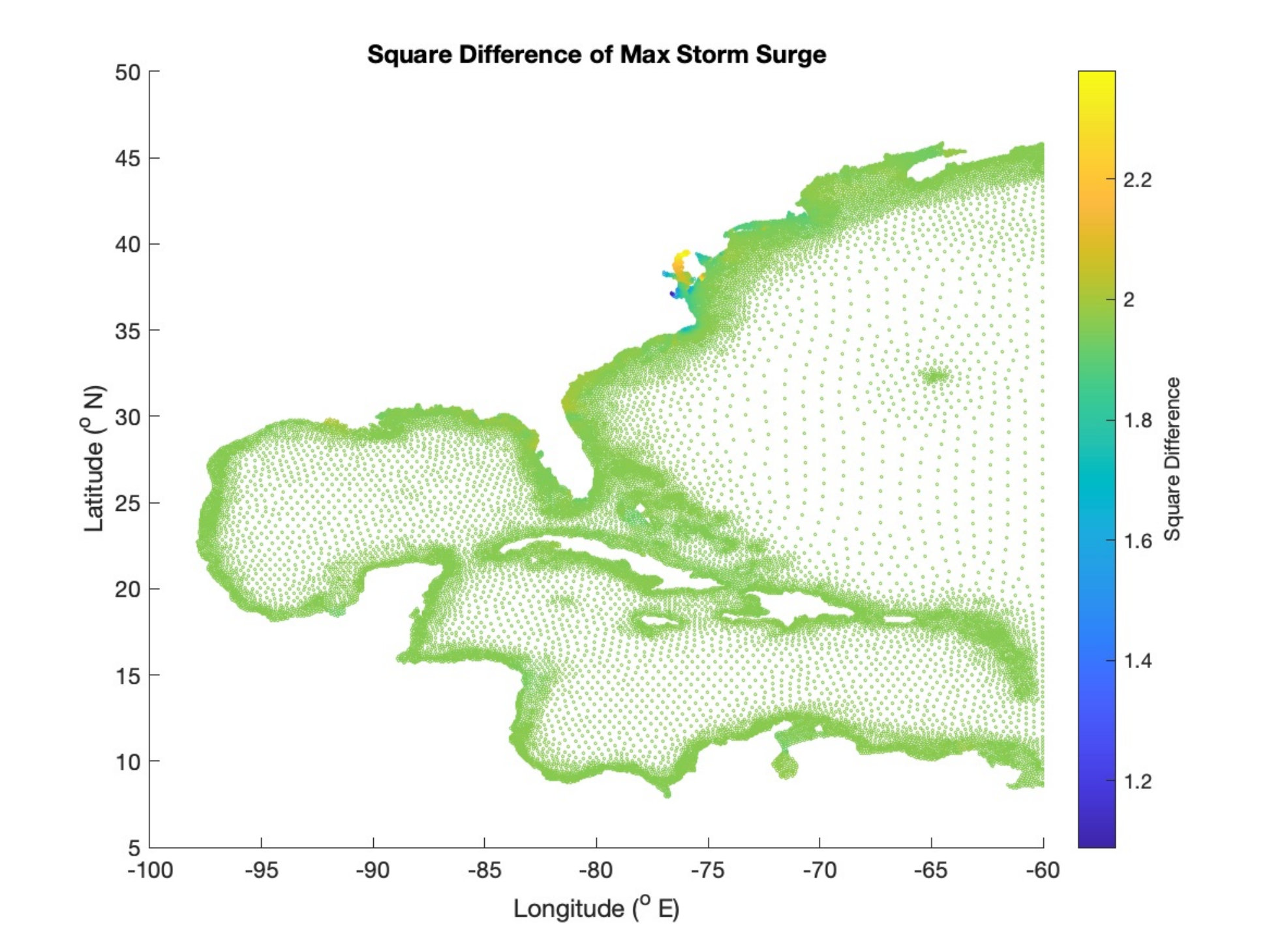}
    \caption{Square Difference from 0 km Storm Surge to 25 km Storm Surge for Hurricane Isabel.}
    \label{fig:stormsurge}
\end{wrapfigure}

\begin{figure*}[htp]
   \subfloat[]{\label{fig:isabelalltime}
      \includegraphics[width=.46\textwidth]{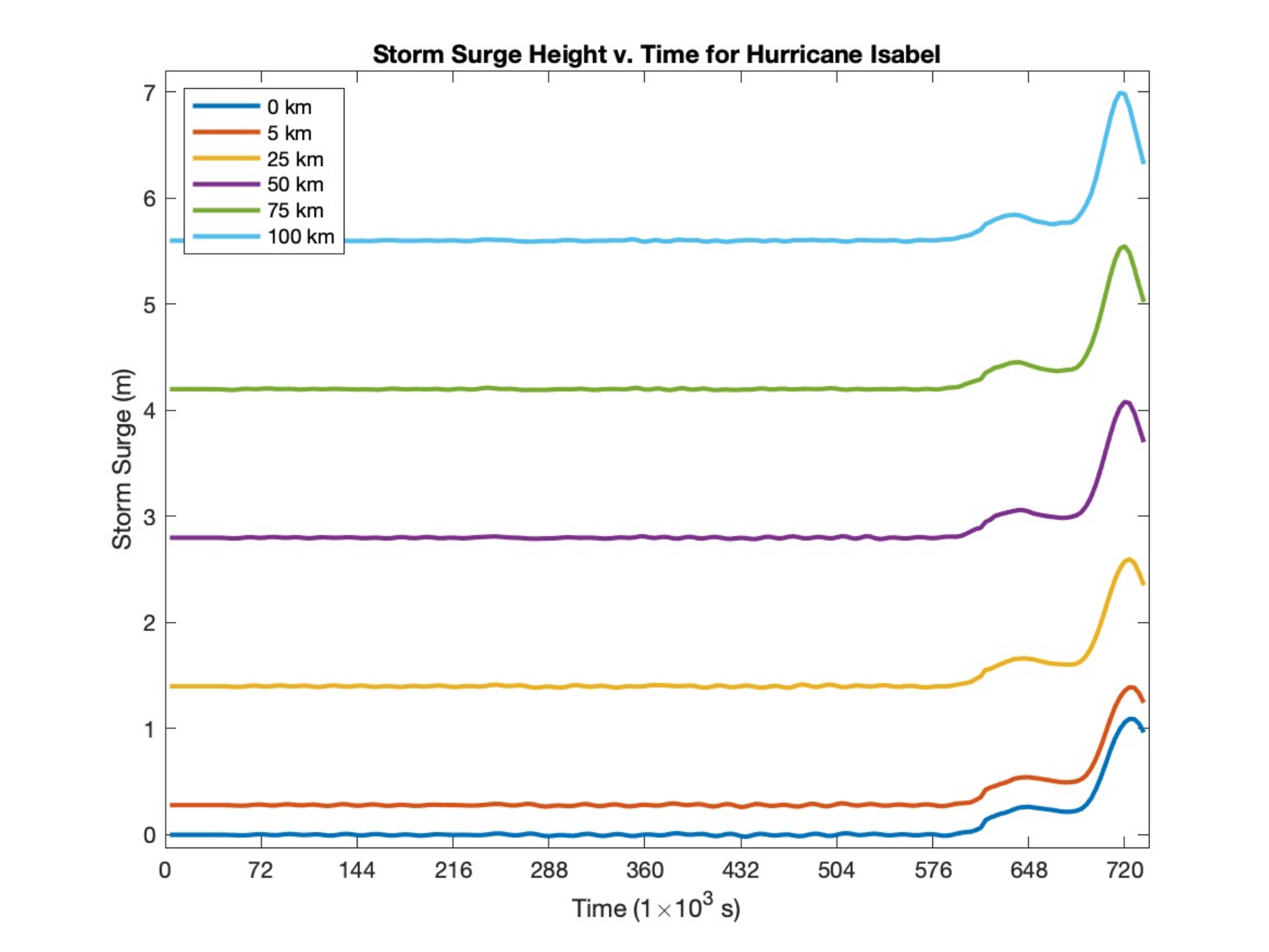}}
~
   \subfloat[]{\label{fig:isabelselect}
      \includegraphics[width=.46\textwidth]{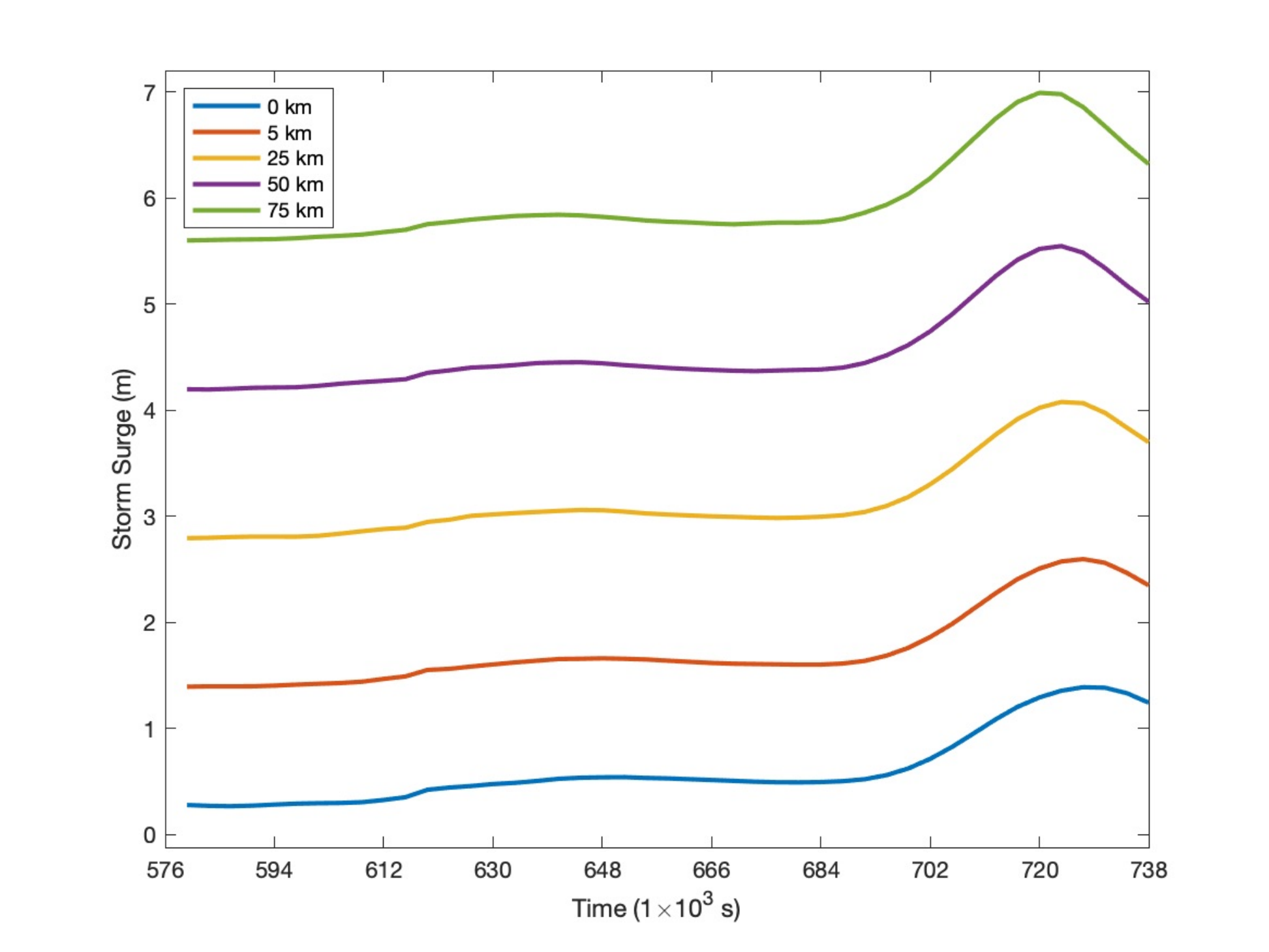}}

   \caption{Storm Surge level over time for our location in the Chesapeake Bay (a) and a closer look at the times when Storm Surge began to take effect (b).}
\end{figure*}

\indent Therefore, we decide to choose a specific location in the Chesapeake Bay Region in order to observe the affect of sea level rise on storm surge levels. The location is 39.132082$^o$ N -76.331032$^o$ E. We see in Figure \ref{fig:isabelselect} that this location is not affected until around time 576, which is likely when Hurricane Isabel approached the region. As before, we see that the model has similar behavior across all sea level rise simulations, but instead it appears that peak storm surge has greater magnitude during time of highest storm surge. 

\section{Discussion}

\subsection{Data Assimilation impact on glacier studies}
\indent Using mean square difference calculations, we decide the best observation scheme and ensemble sizes for our chosen ``true" and ``inaccurate" parameters and initial conditions (Figure \ref{ourParameters}). Indeed, we see in our overall model results that utilizing EnKF improves the output of the model, moving the overall results closer to the ``truth." While this model is simplified, this suggests that using a similar data assimilation technique on a more complex model could help determine the quantity of measurements needed to accurately model ice sheets, and specifically that the increased costs from more observations does not increase accuracy of the model. This is important as it helps glacier modelers to determine the frequency of measurements that researchers should gather in order to get accurate model output. We found that in Section \ref{observationScheme} there should be observations every 19 years before 1900 and yearly after 1950 in order best run the model with minimal error.\\
\indent Our results show that once a threshold is reached for assimilation, the change in mean square difference is small for all further assimilation, so performing additional assimilations does not significantly improve the quality of the models analysis. For data assimilated in the post-satellite era, this threshold is assimilating data once per year. Using a greater number of assimilations is computationally more expensive, making it worthwhile to find an optimal assimilation scheme. More significant is the fact that it can be very expensive and time-intensive for researchers to collect glacier field data. We do not recommend that researchers stop collecting data more frequently than once a year as this data still could provide insight for different areas of research relating to glaciers and the systems they affect. Further, assimilation frequency should be further explored using more complicated models and/or using real observation data. \\
\indent A more complicated model will also require larger ensemble sizes, something that can be computationally expensive to utilize. Therefore, we conjecture that it would be worthwhile to explore a more computationally efficient data assimilation method such as the singular evolutive interpolated Kalman filter (see \cite{hoteit2002} for more details) in more complex models. \\
\indent The results of the sensitivity analysis show that the model itself is quite sensitive to variations in measurements, which is significant as measurements can often have a lot of error in the real world. For one, the knowledge of what range of accuracy in measurement is critical for each parameter may help guide instrumentation development. This underscores the importance of accurate measurements being used in our model, and further the important role that data assimilation may play in readjusting the model runs in accordance with the variations created by these errors.
\subsection{Storm Surge}
\indent We have found that there is a connection between the sea level rise caused by glacial melt and possible increases in storm surge in the future. It is worth noting that while the relationship between these factors appear linear in Figure \ref{fig:stormsurge} and the following figures, these simulations were modeled using a limited number of data points and in very specific cases, so it should not be assumed as a linear relationship. Further, the ADCIRC model physics is based on equations which are non linear, an additional reason not to make that assumption without a more rigorous treatment. This nonlinearity is especially clear when we look at storm surge for a specific grid point as there is both an increase and decrease in magnitude of storm surge for different sea level rise situations and in diffent locations. This shows the need for storm surge simulations to be produced for specific locations in order to see how sea level rise affects storm surge. In these simulations, we have shown that it would be worthwhile to use more complex glacier models to project sea level rise in the future and in turn calculate more accurate possible changes in storm surge. This is important in planning management of possible effects of climate change, especially those due to the increase in glacial melt.\\
\indent A future research direction is exploration of possible increases in storm inundation. While some grid points in this study are over land, we did not make an effort to distinguish these points or to find ways to calculate resulting inundation where such grid points do not exist. It would be worthwhile to inspect changes in coastal inundation with regards to glacial melt driven sea level rise in future studies. \\
\subsection{Study Limitations}
\indent It is worth noting that these calculations were performed using one specific set of parameters, so the findings are likely not universal. Further, when this method is implemented with more complicated glacier models it will still be worthwhile to find this observation scheme, and to determine it with more precision as small differences in observations number may be significant due to non-linearity. \\ 

\section{Conclusions and Conjectures}
In this work, we conducted sensitivity analyses, and investigated and compared the impacts of data assimilation, specifically the Ensemble Kalman Filter (EnKF), on a two-stage ice sheet model. We explore the variables height and length, in the context of glacier melt due to climate change. To the best of the authors’ knowledge, this is one of the only studies in which EnKF techniques are incorporated with ice sheet modeling. The experimental results showed that incorporating observational data with the EnKF technique improves the reliability of glacier height and length forecasts, and we determined that an ensemble size of 10 best minimizes computational cost while also minimizing error. The assimilation results also suggest the time frame in which observational data is most impactful to the output of the model, and that an ice sheet model may perform reasonably well even when only receiving modern observational data once yearly, without any historic data. \\

\indent This work provides stimulus to implement a non-linear data assimilation scheme like the EnKF with a more complicated glacier model and real-world observations. From our ensemble size results (Section \ref{ensembleSize}), we hypothesise that there is a relationship between the number of variables as well as the parameter sensitivity in the two-stage glacier model and the small ideal ensemble size required. When this method is implemented with more complicated glacier models, the ideal ensemble size will likely be larger. We further postulate that a model supported by an observation scheme with a sufficient number of recent observations does not require much pre-satellite era observation data to make good predictions. ``Sufficient" may be as low as one observation per year. However, if the model trajectory is quite far removed from the mean observation trajectory, long-term predictions will generally not be reliable using an EnKF alone.


\section{Acknowledgements}
This work is supported in part by the US National Science Foundation awards DMS-2051019 and DMS-1751636. We would like to thank Dr. Talea Mayo for her mentorship throughout this project and Dr. Alexander Robel for his guidance and ice sheet model. We would also like to thank Dr. Lars Ruthotto and the faculty involved in the Emory University Computational Mathematics for Data Science REU for their invaluable support.

\section{Statements and Declarations}
On behalf of all authors, the corresponding author states that there is no conflict of interest.





\bibliography{manuscript}


\begin{thebibliography}{10}
\ifx \bisbn   \undefined \def \bisbn  #1{ISBN #1}\fi
\ifx \binits  \undefined \def \binits#1{#1}\fi
\ifx \bauthor  \undefined \def \bauthor#1{#1}\fi
\ifx \batitle  \undefined \def \batitle#1{#1}\fi
\ifx \bjtitle  \undefined \def \bjtitle#1{#1}\fi
\ifx \bvolume  \undefined \def \bvolume#1{\textbf{#1}}\fi
\ifx \byear  \undefined \def \byear#1{#1}\fi
\ifx \bissue  \undefined \def \bissue#1{#1}\fi
\ifx \bfpage  \undefined \def \bfpage#1{#1}\fi
\ifx \blpage  \undefined \def \blpage #1{#1}\fi
\ifx \burl  \undefined \def \burl#1{\textsf{#1}}\fi
\ifx \doiurl  \undefined \def \doiurl#1{\url{https://doi.org/#1}}\fi
\ifx \betal  \undefined \def \betal{\textit{et al.}}\fi
\ifx \binstitute  \undefined \def \binstitute#1{#1}\fi
\ifx \binstitutionaled  \undefined \def \binstitutionaled#1{#1}\fi
\ifx \bctitle  \undefined \def \bctitle#1{#1}\fi
\ifx \beditor  \undefined \def \beditor#1{#1}\fi
\ifx \bpublisher  \undefined \def \bpublisher#1{#1}\fi
\ifx \bbtitle  \undefined \def \bbtitle#1{#1}\fi
\ifx \bedition  \undefined \def \bedition#1{#1}\fi
\ifx \bseriesno  \undefined \def \bseriesno#1{#1}\fi
\ifx \blocation  \undefined \def \blocation#1{#1}\fi
\ifx \bsertitle  \undefined \def \bsertitle#1{#1}\fi
\ifx \bsnm \undefined \def \bsnm#1{#1}\fi
\ifx \bsuffix \undefined \def \bsuffix#1{#1}\fi
\ifx \bparticle \undefined \def \bparticle#1{#1}\fi
\ifx \barticle \undefined \def \barticle#1{#1}\fi
\bibcommenthead
\ifx \bconfdate \undefined \def \bconfdate #1{#1}\fi
\ifx \botherref \undefined \def \botherref #1{#1}\fi
\ifx \url \undefined \def \url#1{\textsf{#1}}\fi
\ifx \bchapter \undefined \def \bchapter#1{#1}\fi
\ifx \bbook \undefined \def \bbook#1{#1}\fi
\ifx \bcomment \undefined \def \bcomment#1{#1}\fi
\ifx \oauthor \undefined \def \oauthor#1{#1}\fi
\ifx \citeauthoryear \undefined \def \citeauthoryear#1{#1}\fi
\ifx \endbibitem  \undefined \def \endbibitem {}\fi
\ifx \bconflocation  \undefined \def \bconflocation#1{#1}\fi
\ifx \arxivurl  \undefined \def \arxivurl#1{\textsf{#1}}\fi
\csname PreBibitemsHook\endcsname

\bibitem[\protect\citeauthoryear{Camelo et~al.}{2020}]{camelo2020projected}
\begin{botherref}
\oauthor{\bsnm{Camelo}, \binits{J.}},
\oauthor{\bsnm{Mayo}, \binits{T.L.}},
\oauthor{\bsnm{Gutmann}, \binits{E.D.}}:
Projected climate change impacts on hurricane storm surge inundation in the
  coastal united states.
Frontiers in Built Environment
(2020)
\end{botherref}
\endbibitem

\bibitem[\protect\citeauthoryear{Robel et~al.}{2019}]{robel2019marine}
\begin{barticle}
\bauthor{\bsnm{Robel}, \binits{A.}},
\bauthor{\bsnm{Seroussi}, \binits{H.}},
\bauthor{\bsnm{Roe}, \binits{G.H.}}:
\batitle{Marine ice sheet instability amplifies and skews uncertainty in
  projections of future sea-level rise}.
\bjtitle{PNAS}
\bvolume{116}(\bissue{30}),
\bfpage{14887}--\blpage{14892}
(\byear{2019})
\end{barticle}
\endbibitem

\bibitem[\protect\citeauthoryear{Robel}{2015}]{robel2015long}
\begin{barticle}
\bauthor{\bsnm{Robel}, \binits{A.}}:
\batitle{The long future of antarctic melting}.
\bjtitle{Nature}
\bvolume{526}(\bissue{7573}),
\bfpage{327}--\blpage{328}
(\byear{2015})
\end{barticle}
\endbibitem

\bibitem[\protect\citeauthoryear{ECMWF}{2022}]{dataassim}
\begin{botherref}
\oauthor{\bsnm{ECMWF}}:
Data assimilation
(2022)
\end{botherref}
\endbibitem

\bibitem[\protect\citeauthoryear{Barker et~al.}{2012}]{barker2012weather}
\begin{botherref}
\oauthor{\bsnm{Barker}, \binits{D.}},
\oauthor{\bsnm{Huang}, \binits{X.-Y.}},
\oauthor{\bsnm{Liu}, \binits{Z.}},
\oauthor{\bsnm{Auligné}, \binits{T.}},
\oauthor{\bsnm{Zhang}, \binits{X.}},
\oauthor{\bsnm{Rugg}, \binits{S.}},
\oauthor{\bsnm{Ajjaji}, \binits{R.}},
\oauthor{\bsnm{Bourgeois}, \binits{A.}},
\oauthor{\bsnm{Bray}, \binits{J.}},
\oauthor{\bsnm{Chen}, \binits{Y.}},
\oauthor{\bsnm{Demirtas}, \binits{M.}},
\oauthor{\bsnm{Guo}, \binits{Y.-R.}},
\oauthor{\bsnm{Henderson}, \binits{T.}},
\oauthor{\bsnm{Huang}, \binits{W.}},
\oauthor{\bsnm{Lin}, \binits{H.-C.}},
\oauthor{\bsnm{Michalakes}, \binits{J.}},
\oauthor{\bsnm{Rizvi}, \binits{S.}},
\oauthor{\bsnm{Zhang}, \binits{X.}}:
The weather research and forecasting model's community variational/ensemble
  data assimilation system: Wrfda
\textbf{93}(6),
831--843
(2012)
\end{botherref}
\endbibitem

\bibitem[\protect\citeauthoryear{Robel et~al.}{2018}]{robel2018response}
\begin{barticle}
\bauthor{\bsnm{Robel}, \binits{A.A.}},
\bauthor{\bsnm{Roe}, \binits{G.H.}},
\bauthor{\bsnm{Haseloff}, \binits{M.}}:
\batitle{Response of marine-terminating glaciers to forcing: time scales,
  sensitivities, instabilities, and stochastic dynamics}.
\bjtitle{Journal of Geophysical Research: Earth Surface}
\bvolume{123}(\bissue{9}),
\bfpage{2205}--\blpage{2227}
(\byear{2018})
\end{barticle}
\endbibitem

\bibitem[\protect\citeauthoryear{Asch et~al.}{2009}]{asch2016fundamentals}
\begin{bbook}
\bauthor{\bsnm{Asch}, \binits{M.}},
\bauthor{\bsnm{Bocquet}, \binits{M.}},
\bauthor{\bsnm{Nodet}, \binits{M.}}:
\bbtitle{Fundamentals of Algorithms Data Assimilation}.
\bpublisher{SIAM}, \blocation{???}
(\byear{2009})
\end{bbook}
\endbibitem

\bibitem[\protect\citeauthoryear{Loucks et~al.}{2005}]{modelanalysis}
\begin{bchapter}
\bauthor{\bsnm{Loucks}, \binits{D.P.}},
\bauthor{\bsnm{Beek}, \binits{E.v.}},
\bauthor{\bsnm{Stedinger}, \binits{J.R.}},
\bauthor{\bsnm{Dijkman}, \binits{J.P.M.}},
\bauthor{\bsnm{Villars}, \binits{M.T.}}:
\bctitle{Model Sensitivity and Uncertainty Analysis}.
\bbtitle{Water Resources Systems Planning and Management},
pp. \bfpage{255}--\blpage{286}.
\bpublisher{UNESCO}, \blocation{???}
(\byear{2005})
\end{bchapter}
\endbibitem

\bibitem[\protect\citeauthoryear{Bjork et~al.}{2015}]{bjork2015brief}
\begin{barticle}
\bauthor{\bsnm{Bjork}, \binits{A.A.}},
\bauthor{\bsnm{Kruse}, \binits{L.M.}},
\bauthor{\bsnm{Michaelsen}, \binits{P.B.}}:
\batitle{Brief communication: Getting greenland's glaciers right – a new
  dataset of all official greenlandic glacier names}.
\bjtitle{The Cryosphere Discuss}
\bvolume{9},
\bfpage{1593}--\blpage{1601}
(\byear{2015})
\end{barticle}
\endbibitem

\bibitem[\protect\citeauthoryear{Hotiet et~al.}{2002}]{hoteit2002}
\begin{barticle}
\bauthor{\bsnm{Hotiet}, \binits{I.}},
\bauthor{\bsnm{Pham}, \binits{D.-T.}},
\bauthor{\bsnm{Blum}, \binits{J.}}:
\batitle{A simplified reduced order kalman filtering and application to
  altimetric data assimilation in tropical pacific}.
\bjtitle{Jounral of Marine Systems}
\bvolume{36},
\bfpage{101}--\blpage{127}
(\byear{2002})
\end{barticle}
\endbibitem

\end{thebibliography}


\onecolumn
\appendix
\section{Code}
The full data file with the data we generated, as well as related code and further graphs of data assimilation schemes for this two-stage ice sheet model can be found at 
\url{https://github.com/hakuupi/StormSurge}. \\

The code used to run model and perform assimilation is in the {\fontfamily{pcr}\selectfont {\small \textbf{Juypter\_Notebook\_Code}}} folder; we ran Python in Jupyter Notebook to generate data and run model.\\

Our data is stored in {\fontfamily{pcr}\selectfont {\small \textbf{Experiment\_Data}}} and is saved in csv files. Note that the parameters and initial conditions used are in the name of the file.\\

The MATLAB code used to generate graphics used in our paper is found in {\fontfamily{pcr}\selectfont {\small \textbf{Matlab\_Plot\_Code}}}.

\newpage
\section{Figures}

\begin{figure*}[htp]
   \subfloat[Assimilation Frequencies Before 1900 in the pre-satellite era. The assimilation schemes with under 5 percent in the assimilation time period are every 19 years or more frequent. The scheme that assimilates data every 19 years is boxed. The red vertical line marks the year 1900, the end of data assimilation in this experiment.]{\label{fig:oldDataGraphs}
      \includegraphics[width=\textwidth]{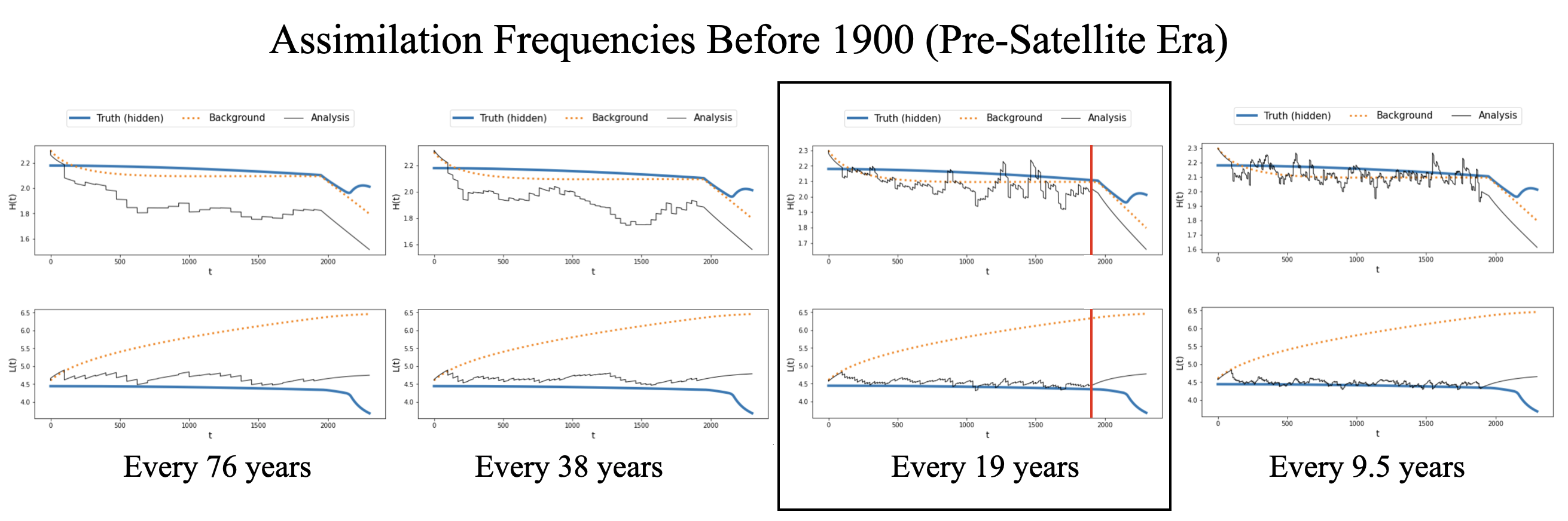}}

   \subfloat[Assimilation Frequencies After 1950 in the post-satellite era. The assimilation schemes with under 5 percent in the assimilation time period are yearly or more frequent. The yearly scheme is boxed. The red vertical line marks the year 1950, the start of data assimilation in this experiment.]{\label{fig:newDataGraphs}
      \includegraphics[width=\textwidth]{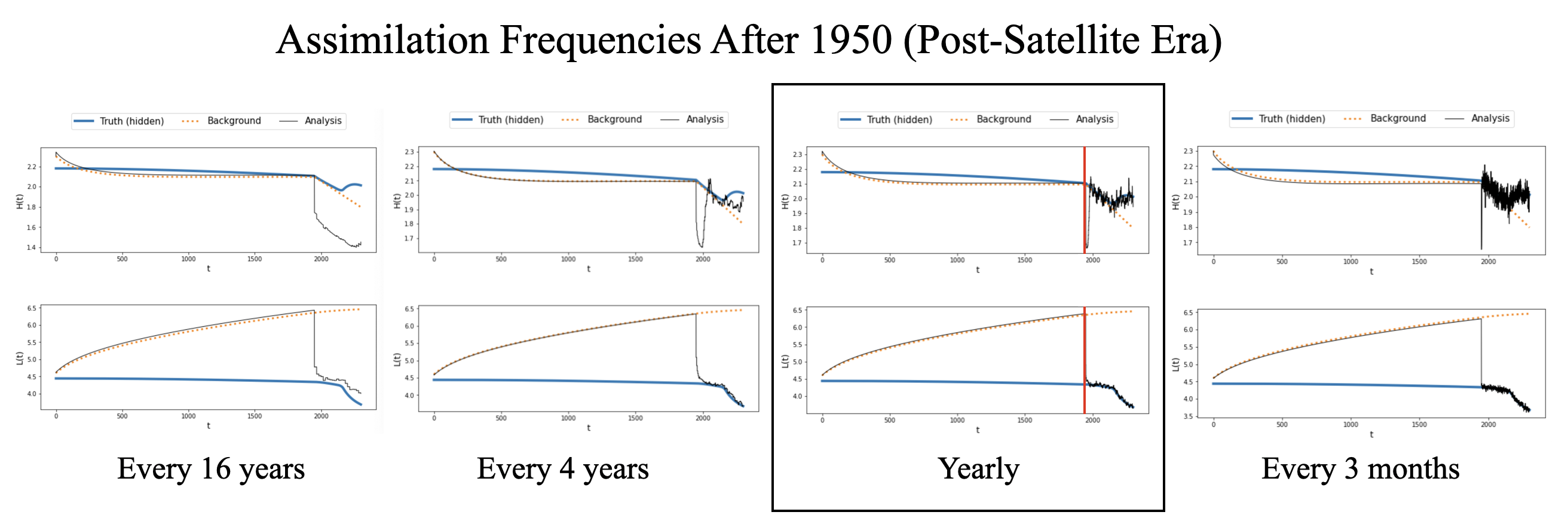}}
   \caption{Assimilation dates and frequencies}\label{fig: assimilationDates}
\end{figure*}

\begin{figure*}[htp]
   \subfloat[Observations and Truth Simulation for H]{
      \includegraphics[width=.48\textwidth]{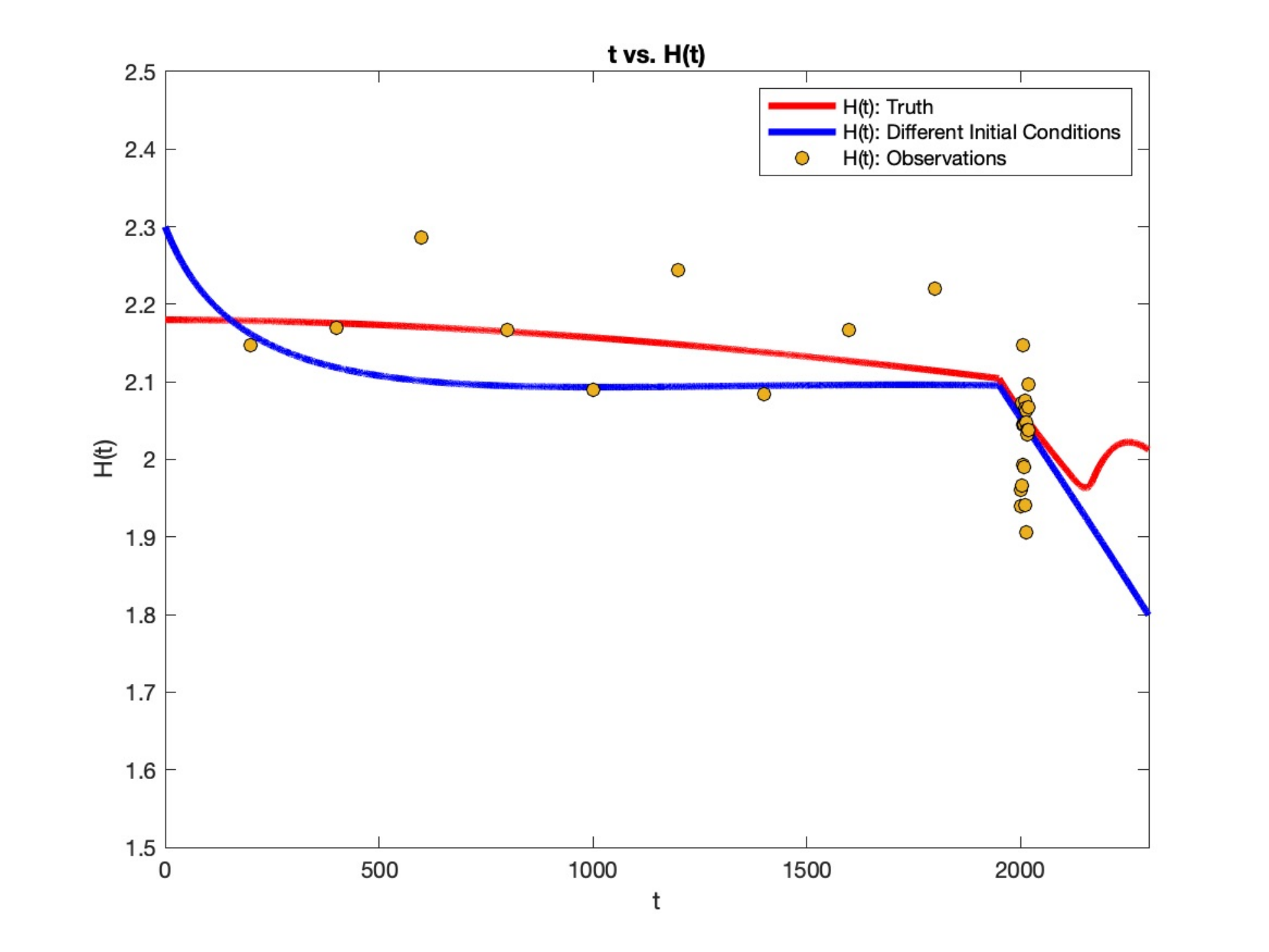}}
~
   \subfloat[Observations and Truth Simulation for L]{
      \includegraphics[width=.48\textwidth]{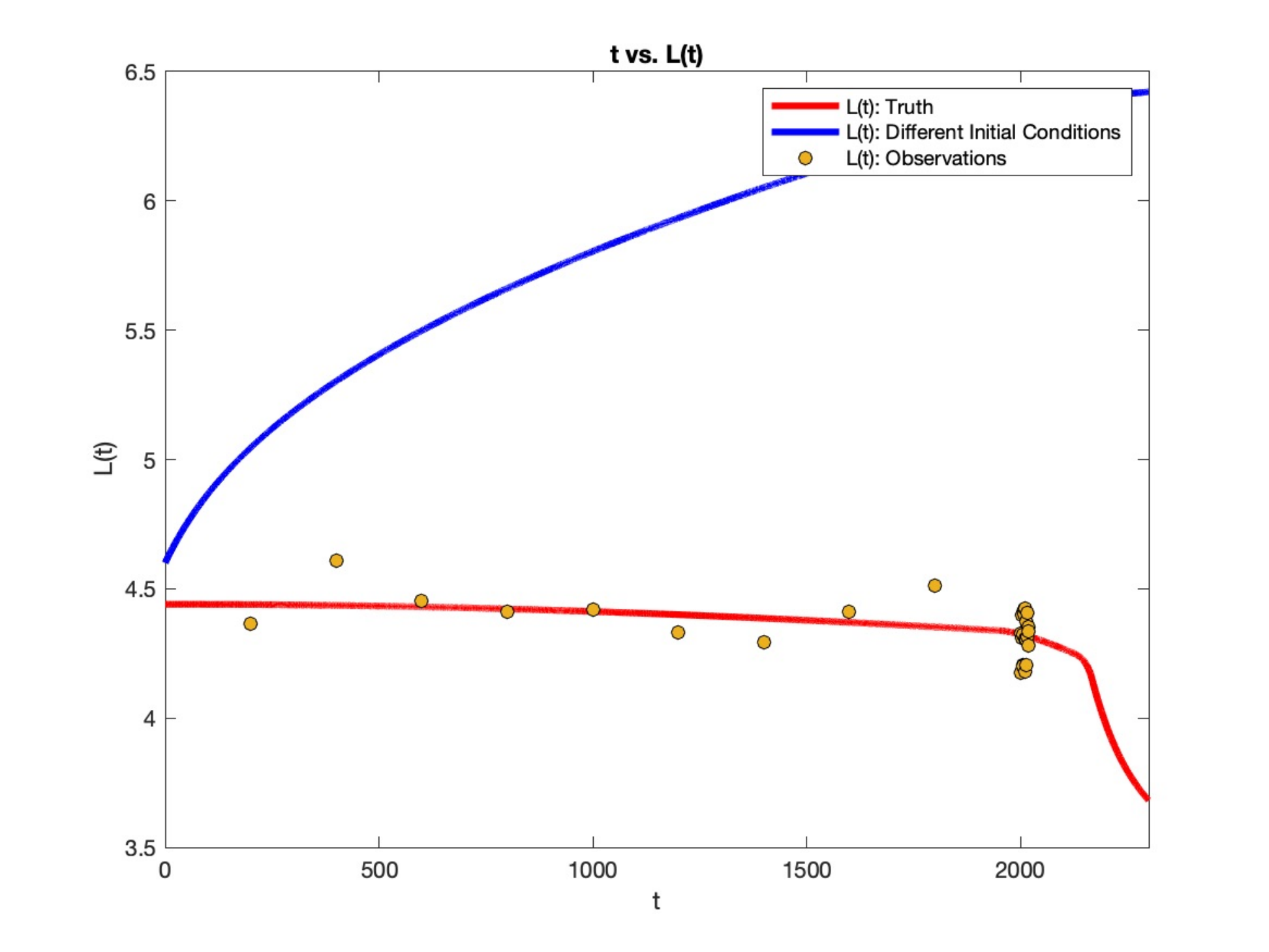}}

   \subfloat[Analysis, Ensemble, Observations and Truth Simulation for H]{
      \includegraphics[width=.48\textwidth]{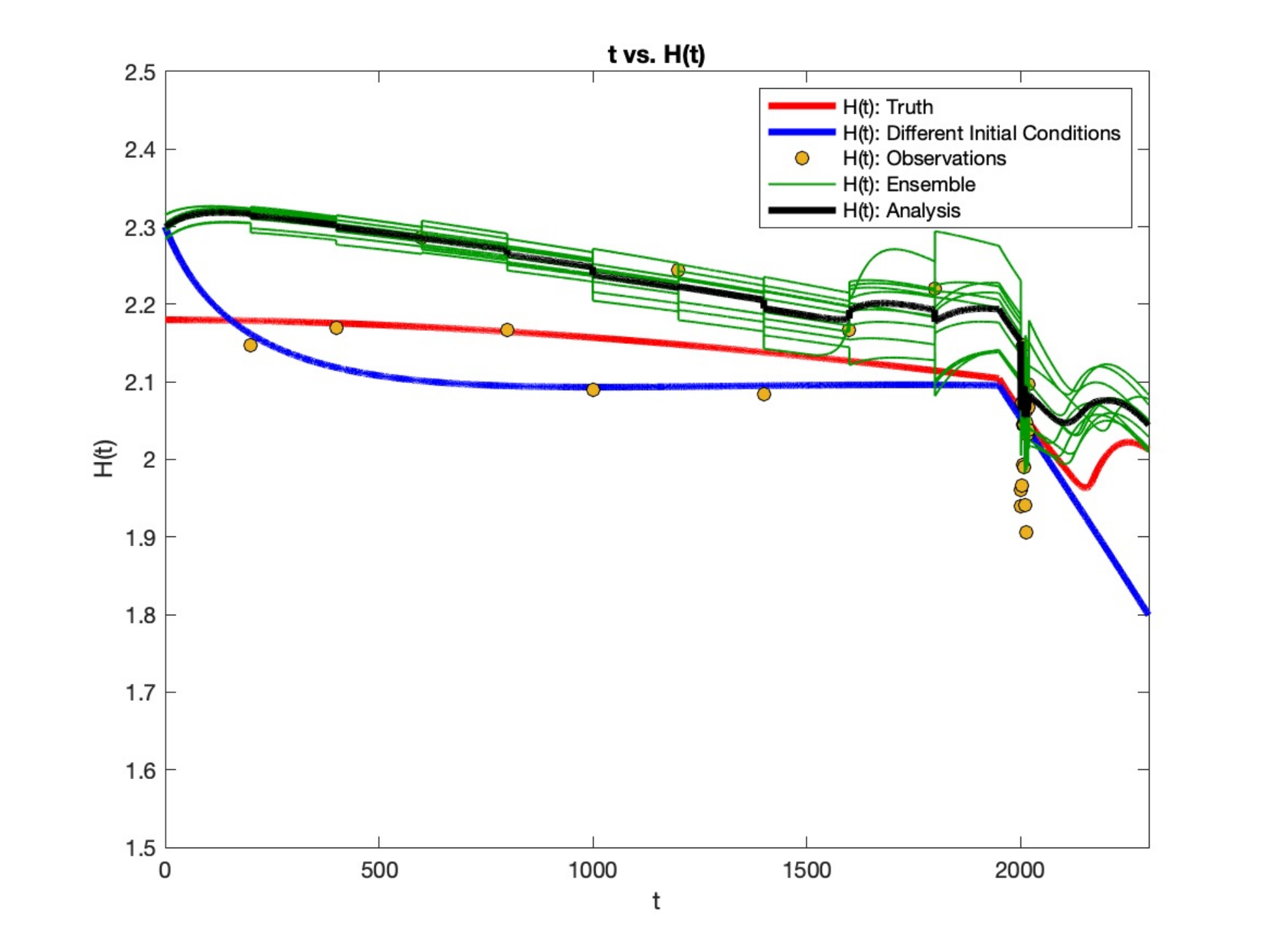}}
~
   \subfloat[Analysis, Ensemble, Observations and Truth Simulation for L]{
      \includegraphics[width=.48\textwidth]{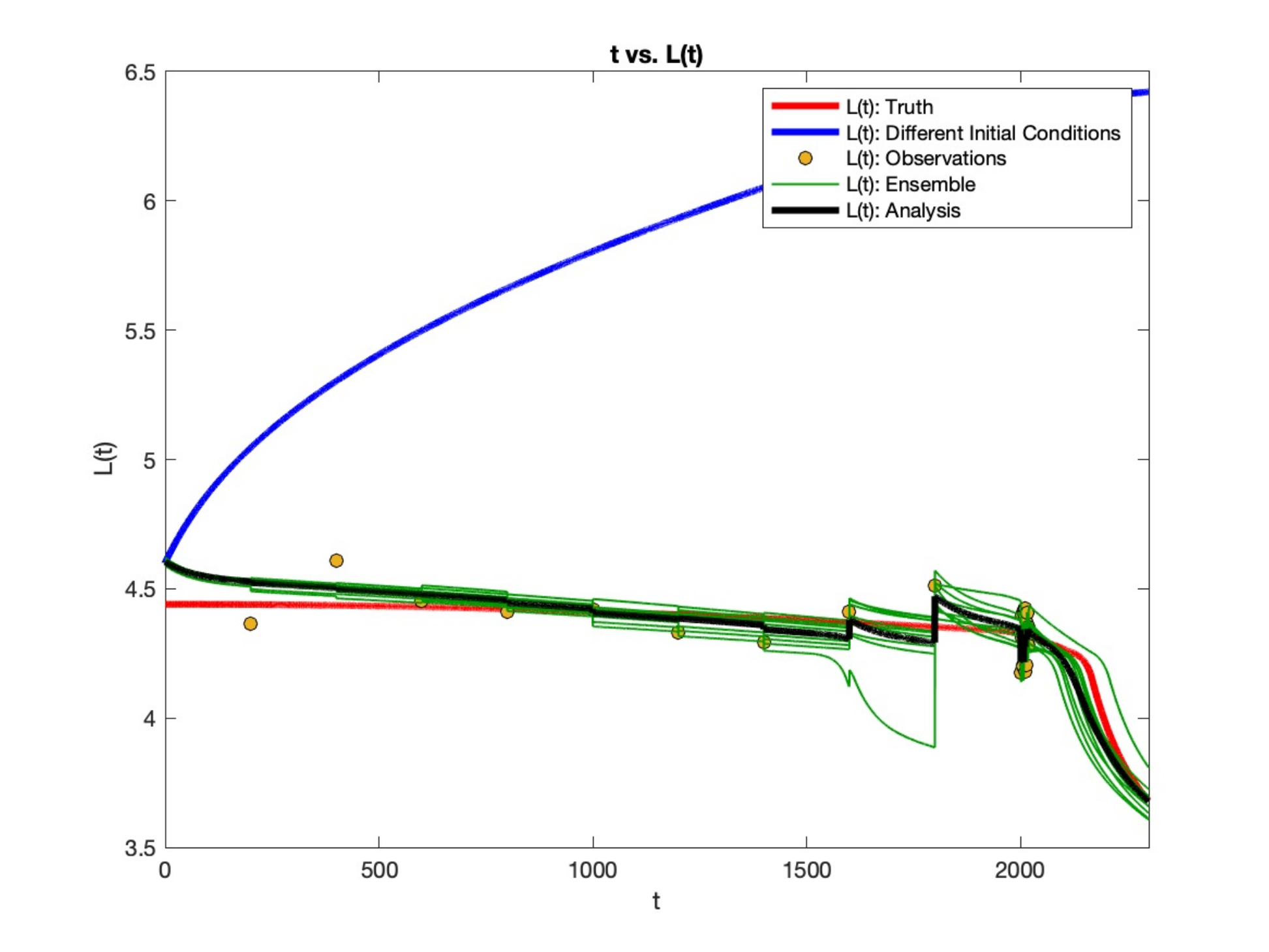}}

     \subfloat[Analysis, Observations and Truth Simulation for H]{
      \includegraphics[width=.48\textwidth]{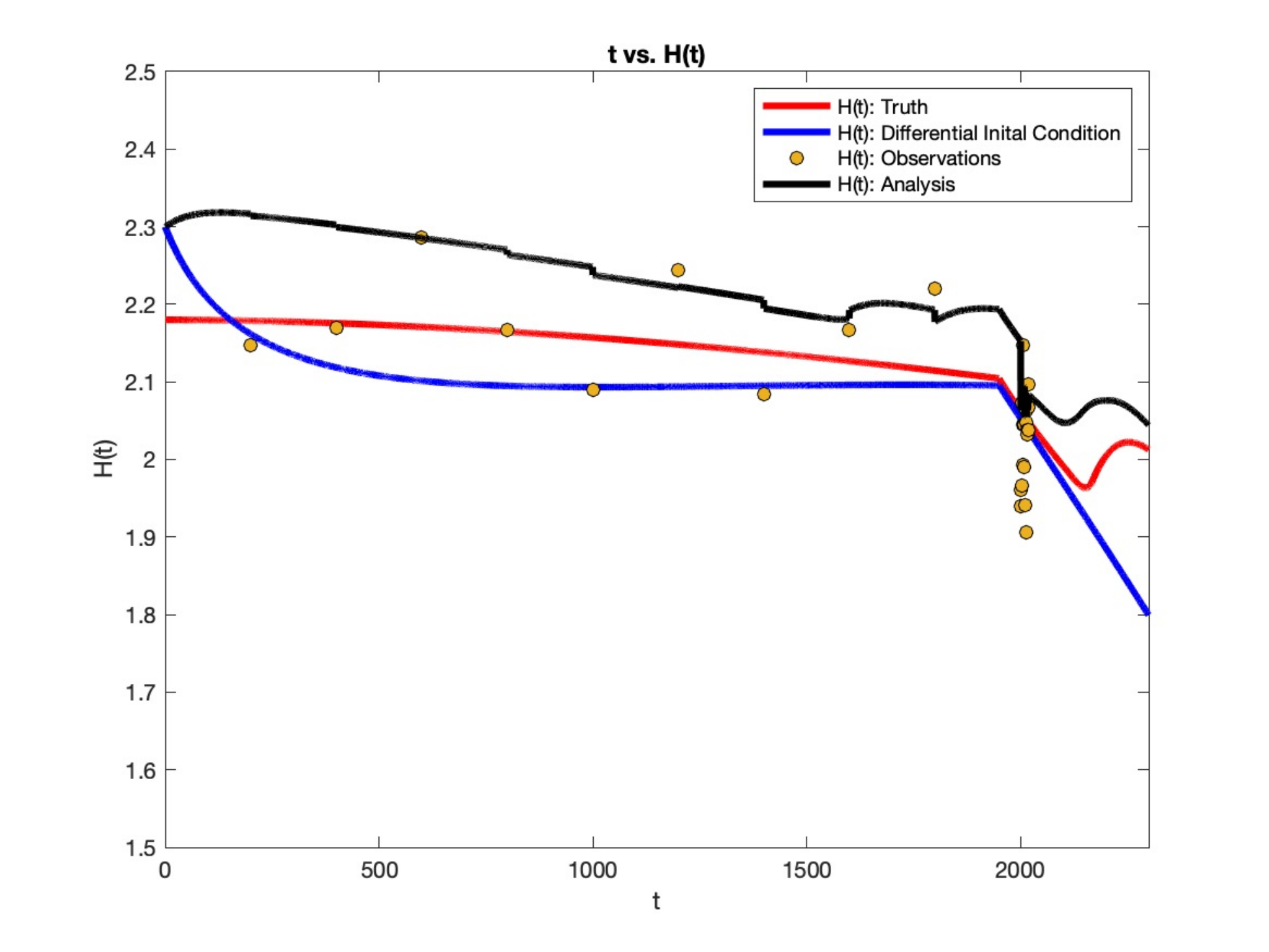}}
~
   \subfloat[Analysis, Observations and Truth Simulation for L]{
      \includegraphics[width=.48\textwidth]{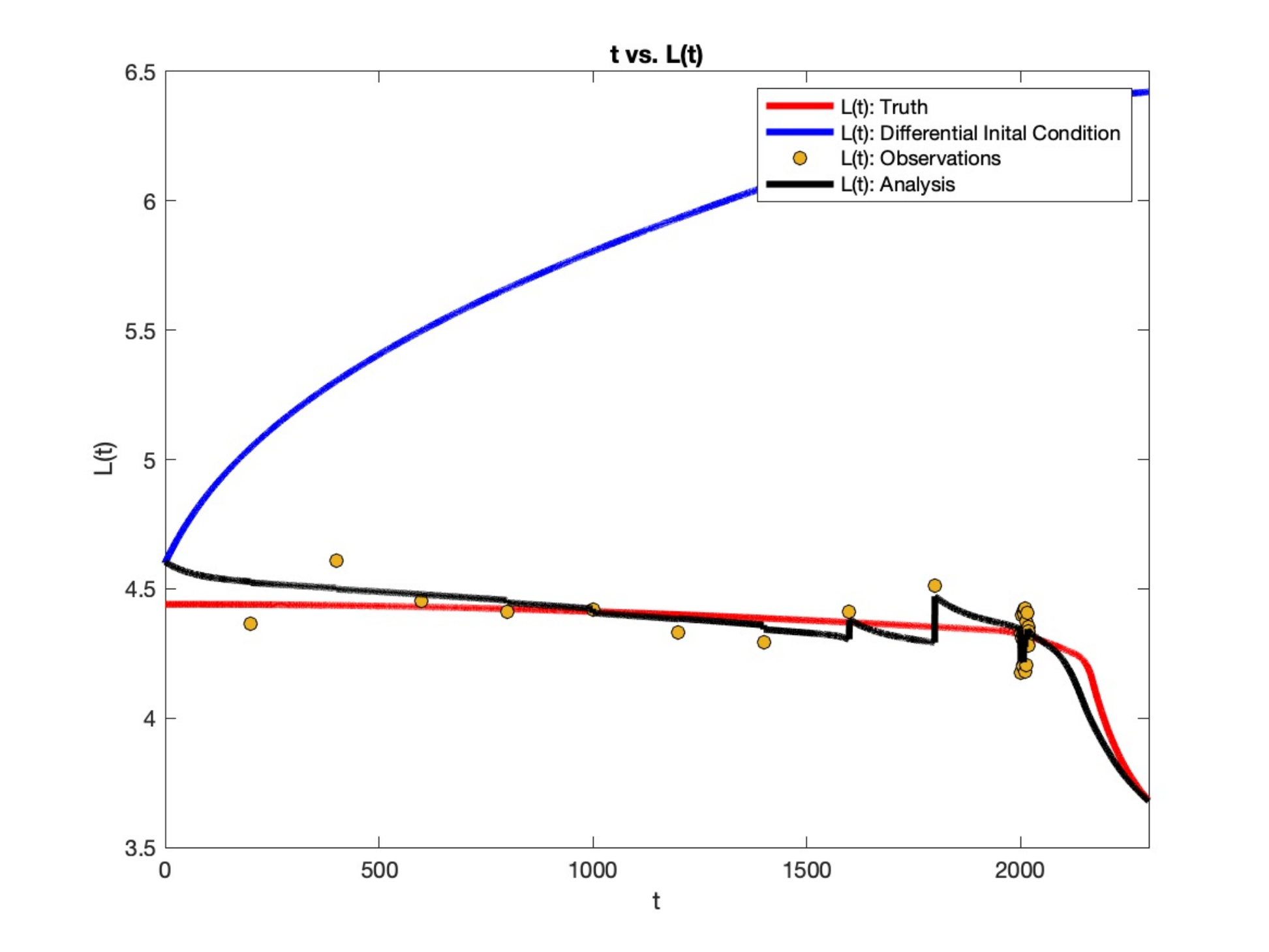}}    
  
   \caption{Model runs with worse observation and ensemble Scheme with mean square difference for H of 0.00982 and mean square difference for L of 0.00941. Note: mean square difference is rounded to 3 significant digits.}
\end{figure*}

\begin{figure*}[htp]
   \subfloat[Observations and Truth Simulation for H]{
      \includegraphics[width=.48\textwidth]{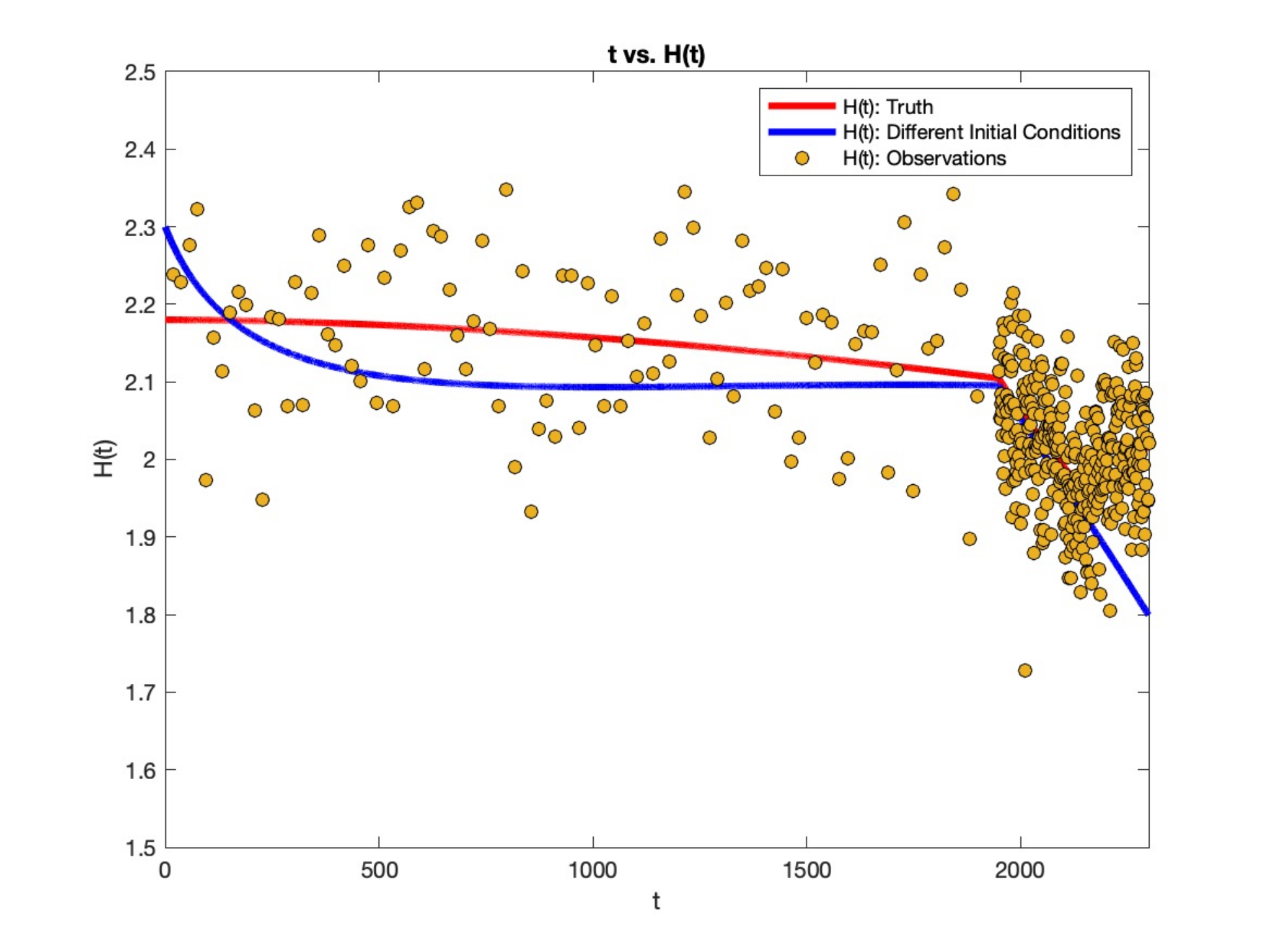}}
      ~
    \subfloat[Observations and Truth Simulation for L]{
      \includegraphics[width=.48\textwidth]{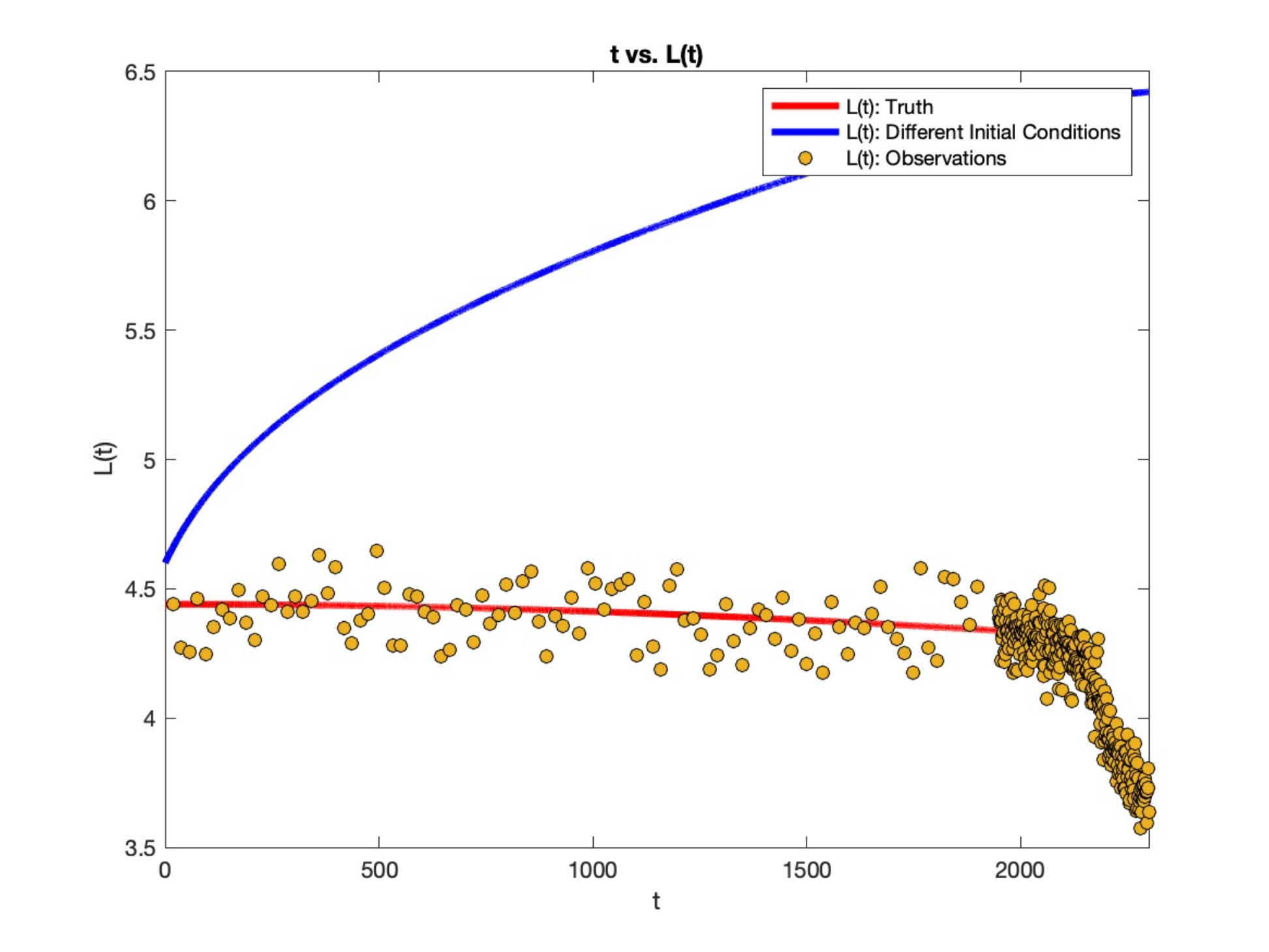}}
      
    \subfloat[Analysis, Ensemble, Observations and Truth Simulation for H]{
      \includegraphics[width=.48\textwidth]{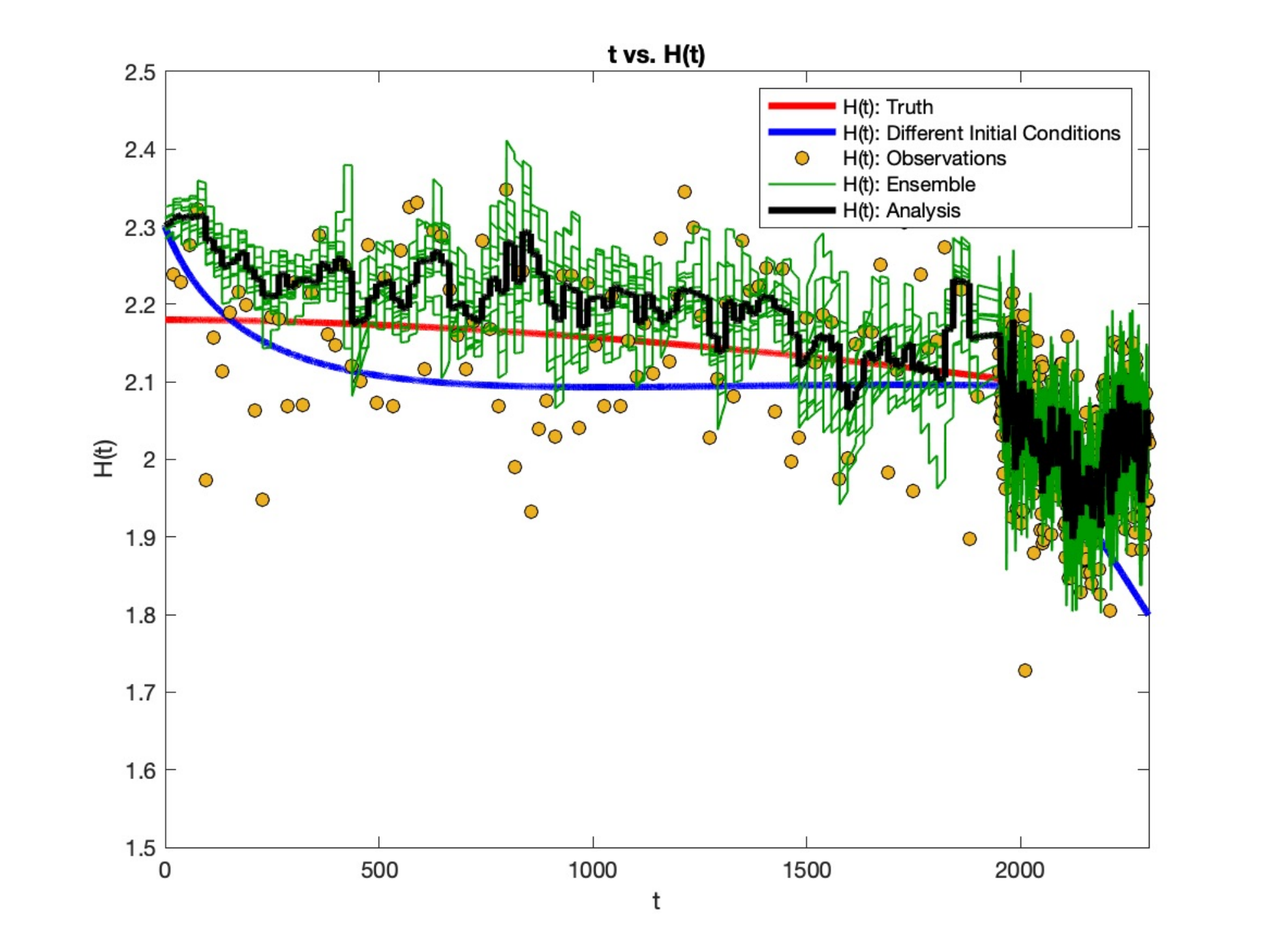}}
      ~
    \subfloat[Analysis, Ensemble, Observations and Truth Simulation for L]{
      \includegraphics[width=.48\textwidth]{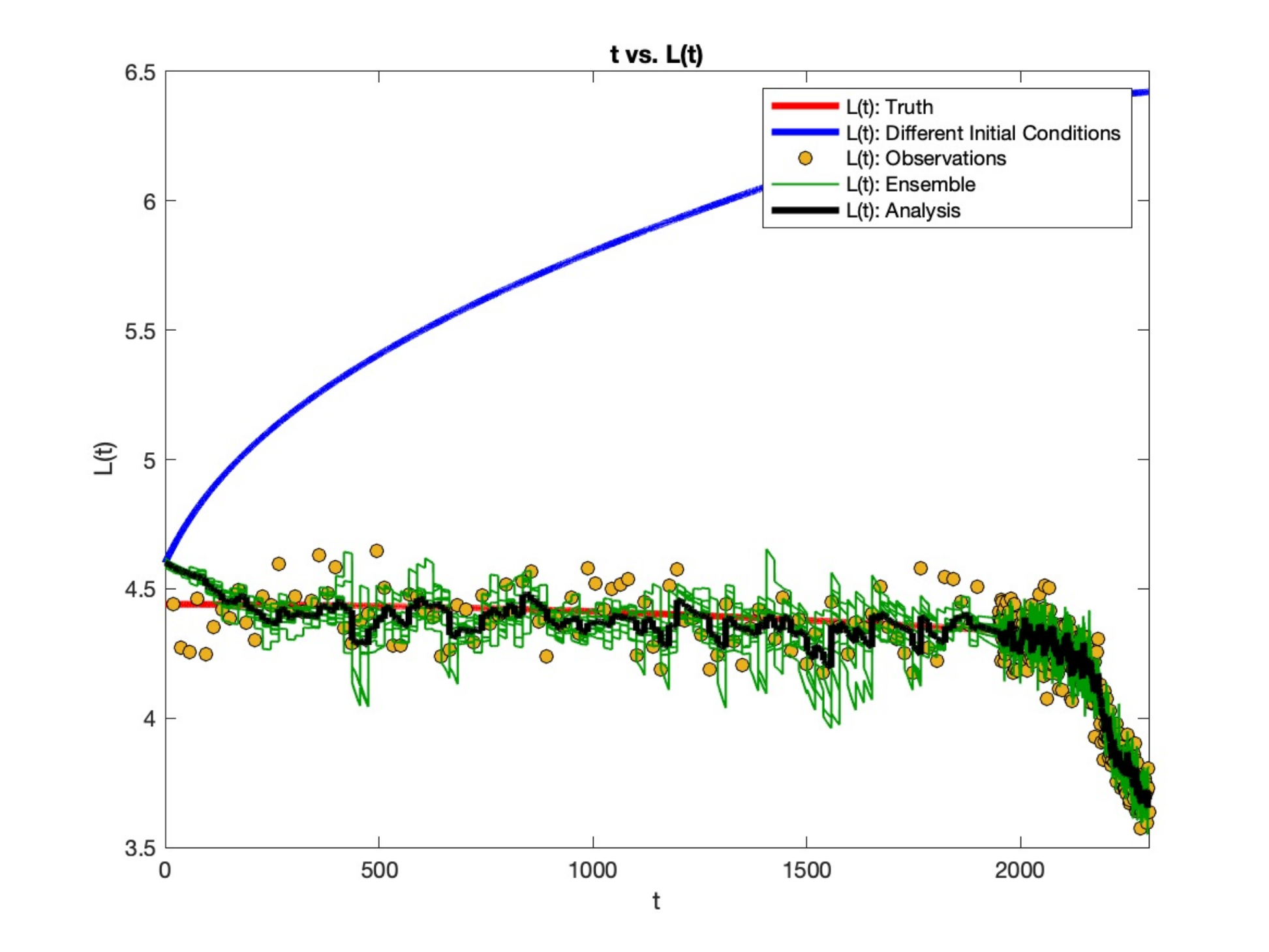}}
      
    \subfloat[Analysis, Observations and Truth Simulation for H]{
      \includegraphics[width=.48\textwidth]{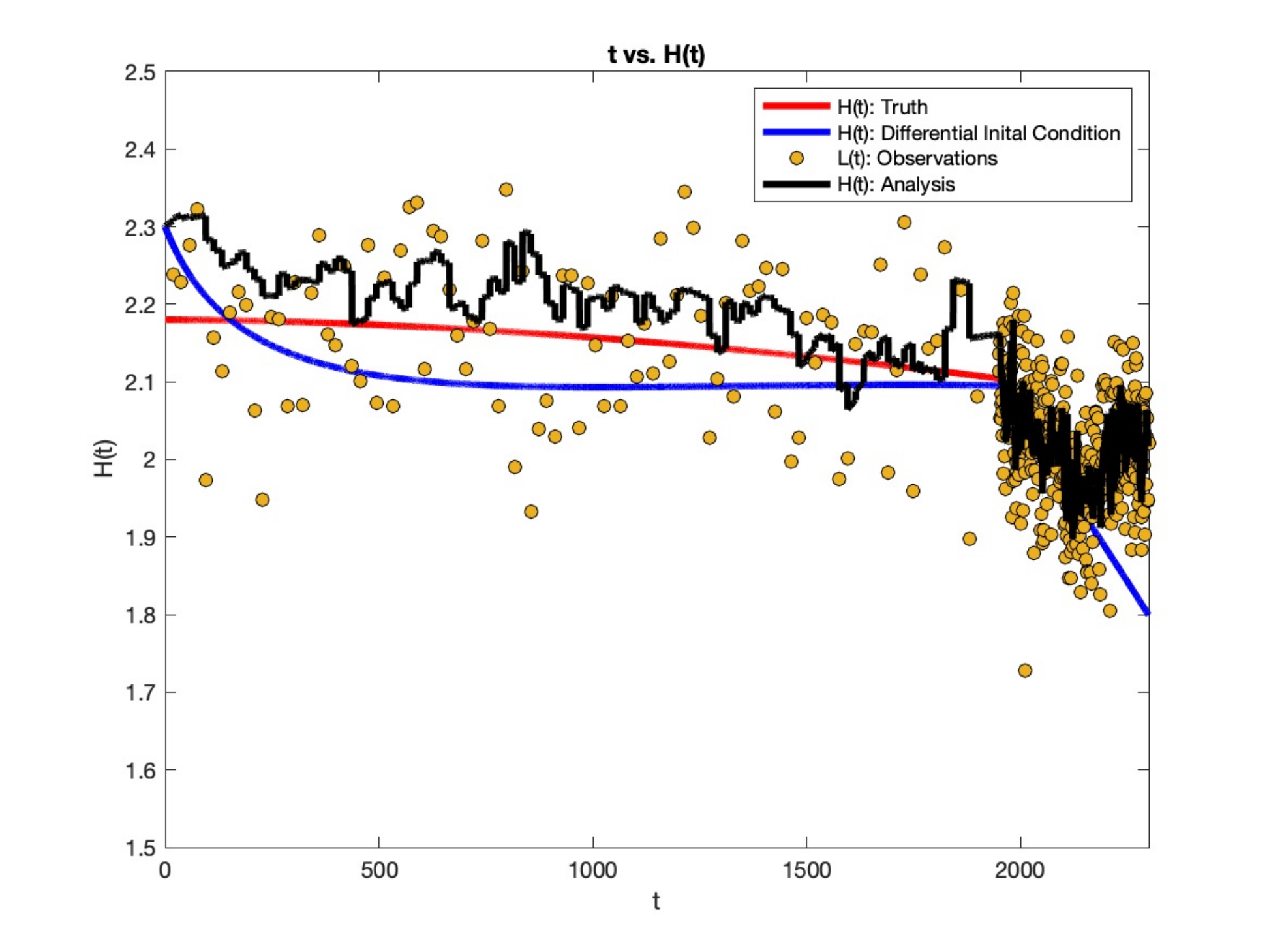}}
      ~
    \subfloat[Analysis, Observations and Truth Simulation for L]{
      \includegraphics[width=.48\textwidth]{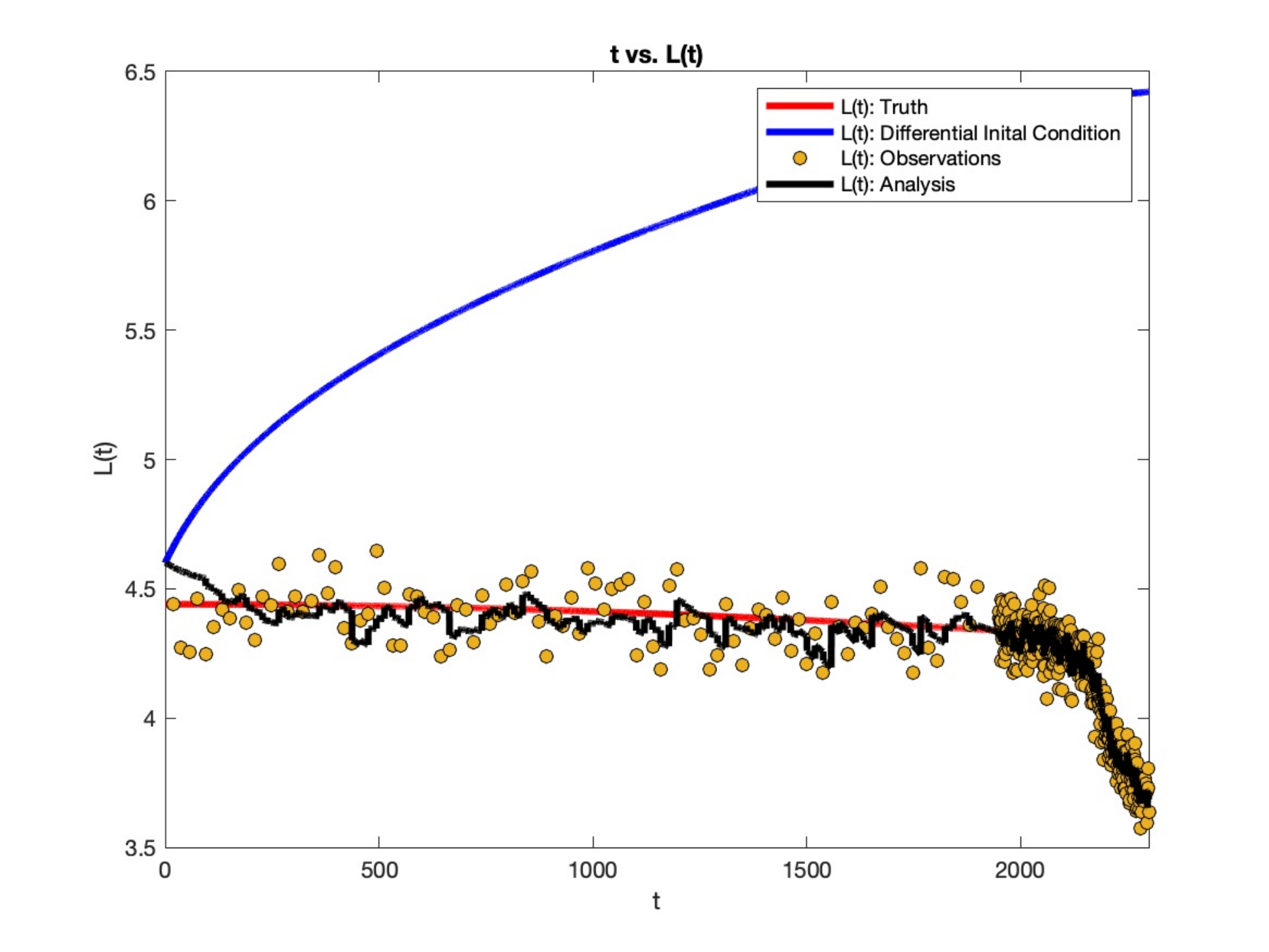}}
      
    \caption{Model Runs with best observation and ensemble scheme, with mean square difference for H of 0.00353 and mean square difference for L of 0.00328.  Note: mean square difference is rounded to 3 significant digits.}
\end{figure*}

\end{document}